\documentclass[preprint,12pt]{elsarticle}

\usepackage{amssymb}
\usepackage{amsmath} 
\usepackage{graphicx}
\usepackage{hyperref}
\usepackage{booktabs}
\usepackage{listings}
\usepackage{xcolor}
\usepackage{ragged2e}
\usepackage{float}
\usepackage{algorithm} 
\usepackage{algpseudocode} 
\usepackage{caption}
\usepackage{setspace}
\usepackage[utf8]{inputenc}
\usepackage{lineno} 

\journal{}

\begin{document}

\begin{frontmatter}


\author[inst1]{Mingjiao Yan}
\author[inst2]{Yang Yang}
\author[inst2]{Zongliang Zhang}
\author[inst1]{Dengmiao Hao}
\author[inst1]{Chao Su}
\author[inst3]{Qingsong Duan}
\affiliation[inst1]{organization={College of Water Conservancy and Hydropower Engineering, Hohai University},
            city={Nanjing},
            postcode={210098}, 
            state={Jiangsu},
            country={China}}

\affiliation[inst2]{organization={PowerChina Kunming Engineering Corporation Limited},
            city={Kunming},
            postode={650051}, 
            state={Yunnan}, 
            country={China}}
\affiliation[inst3]{organization={College of Water Conservancy, Yunnan Agricultural University},
            city={Kunming},
            postcode={650201}, 
            state={Yunnan},
            country={China}}

\title{Three dimensional seepage analysis using a polyhedral scaled boundary finite element method}

\begin{abstract}
This work presents a polyhedral scaled boundary finite element method (PSBFEM) for three dimensional seepage analysis. We first derive the scaled boundary formulation for 3D seepage problems, and subsequently incorporate Wachspress shape functions to construct shape functions over arbitrary polygonal elements, thereby establishing the foundation of the proposed polyhedral SBFEM. The method combines the semi-analytical nature of the SBFEM with the geometric flexibility of polyhedral and octree meshes, making it well-suited for complex seepage simulations. The PSBFEM is implemented within the ABAQUS UEL framework to facilitate steady-state, transient, and free-surface seepage analyses. A series of numerical examples are conducted to verify the accuracy, efficiency, and convergence properties of the proposed approach, including benchmark tests and applications with intricate geometries. The results demonstrate that the PSBFEM achieves higher accuracy and faster convergence than conventional FEM, particularly when using hybrid octree meshes with local refinement. This framework provides a robust and efficient computational tool for three-dimensional seepage analysis in geotechnical and hydraulic engineering applications.
\end{abstract}


\begin{highlights}
\item Proposed a semi-analytical 3D SBFEM formulation that improves accuracy and efficiency in seepage modeling
\item Developed a 3D polyhedral SBFEM framework incorporating Wachspress shape functions for complex domain modeling
\item Implemented the PSBFEM in ABAQUS UEL to solve steady-state and transient seepage problems
\item Introduced a fixed-mesh strategy to enable free-surface modeling within the PSBFEM framework
\item The PSBFEM achieves higher accuracy than FEM under mesh refinement.
\end{highlights}

\begin{keyword}
Polyhedral scaled boundary finite element method \sep 
Three-dimensional seepage analysis \sep 
Hybrid octree mesh \sep 
Free surface flow \sep 
Wachspress shape functions \sep 
ABAQUS UEL \sep 
\end{keyword}

\end{frontmatter}


\section{Introduction}
\label{sec:1}
Seepage analysis plays a crucial role in various engineering fields, including hydrology \cite{scudelerExam2017}, geotechnical engineering \cite{zhouGroundwate2023}, and environmental science \cite{shahSeepageLosses2021}, etc. Accurate seepage modeling is essential for groundwater management, soil contamination prediction, and stability assessment of earth structures, providing critical insights for engineering decision-making. For complex seepage problems, analytical and experimental approaches are often impractical or infeasible, making numerical methods a compelling and practical alternative. 

Currently, the finite element method (FEM) is one of the widely used approaches for solving seepage problems \cite{johariStochas2019,abokwiekFinite2022,liSelect2022}. However, conventional FEM suffers from several limitations when applied to large-scale three-dimensional (3D) seepage problems with intricate geometries. First, mesh generation for hexahedral elements is labor-intensive and often requires manual intervention, while tetrahedral elements—though more amenable to automation—frequently lead to reduced numerical accuracy \cite{schneiderLarge2022}. Furthermore, FEM requires discretization in all spatial directions and relies on weak forms, which may compromise accuracy and computational efficiency when dealing with low-order elements or poorly conditioned meshes \cite{zienkiewic2005}.

To overcome the limitations of traditional mesh-based methods in solving partial differential equations (PDEs), researchers have proposed a variety of meshfree approaches, including the element-free Galerkin (EFG) method \cite{hegen1996element}, smoothed particle hydrodynamics (SPH) \cite{bonet2004variational}, reproducing kernel particle method (RKPM) \cite{huang2020rkpm2d}, radial basis function (RBF) method \cite{ghoneim2020smoothed}, peridynamics \cite{shojaei2022hybrid}, physics-informed neural networks (PINNs) \cite{luo2023novel} and others. Despite the flexibility and robustness of meshfree methods in handling large deformations and discontinuities, they often suffer from high computational cost, difficulty in imposing essential boundary conditions, complex implementation, and challenges in numerical integration and stability.

In contrast, researchers have also developed several enhanced mesh-based methods to overcome the limitations of traditional FEM. These include isogeometric analysis (IGA) \cite{cottrell2009isogeometric,gupta2023insight}, which integrates CAD and analysis by using NURBS-based shape functions; Other notable developments include the polygonal/polyhedral smoothed finite element method (PSFEM) \cite{zeng2018smoothed,yan2024fast}, polygonal/polyhedral finite element method (PFEM) \cite{wu2023polygonal,yang2025steady}, and their combinations with high-order or isogeometric formulations. Polygonal and polyhedral finite elements allow greater flexibility in mesh generation due to their ability to accommodate arbitrary shapes, making them well-suited for complex geometries and adaptive refinement. They can also reduce the total number of elements, improving computational efficiency without sacrificing accuracy \cite{yan2025polyhedral}. These features make them promising for next-generation mesh-based methods that bridge geometric modeling and numerical analysis.

However, such polygonal and polyhedral finite element technologies still rely on weak form discretization in each spatial direction, which may degrade computational accuracy. This limitation motivates the development of alternative formulations such as the scaled boundary finite element method (SBFEM), which transforms partial differential equations into ordinary differential equations along the radial direction, enabling semi-analytical solutions \cite{song2018scaled}. The SBFEM has been successfully extended to a variety of physical problems, including linear elasticity \cite{krome2017semi,yangDevelopmentABAQUSUEL2020}, fracture mechanics\cite{zhang2024adaptive,ooi2010hybrid}, heat conduction\cite{li2016novel,yan2025polyhedral}, nonlinear analysis \cite{chen2018efficient,zhang2020nonlocal}, and wave propagation\cite{gravenkamp2017efficient,zhang2024prismatic}, demonstrating its versatility and robustness across different fields. 

Since the SBFEM only discretizes the boundary of the domain, it naturally accommodates polygonal and polyhedral elements, making it particularly well-suited for complex geometries and generalized meshes \cite{ya2021open,ye2021psbfem}. Dai et al. \cite{dai2015fully} developed an automatic remeshing procedure based on the polygonal SBFEM for modeling arbitrary crack propagation. Natarajan et al. \cite{natarajan2017scaled} developed a scaled boundary finite element method for 3D convex polyhedra by discretizing only the surfaces using Wachspress interpolants. Yan et al. \cite{yan2025polyhedral} proposed a polyhedral SBFEM solving three-dimensional heat conduction problems. These developments demonstrate that the polygonal and polyhedral SBFEM combines the flexibility of general meshes with the accuracy and efficiency of semi-analytical formulations.

To enable fast and automatic analysis, the SBFEM has been combined with octree meshes, allowing complex structures to be efficiently discretized into polyhedral elements.  Chen et al. \cite{chen2018efficient} proposed a nonlinear octree SBFEM with mean-value interpolation and internal tetrahedral integration, enabling fast automatic analysis of complex geotechnical structures. Saputra et al. \cite{saputra2020three} employed octree-based SBFEM with transition elements to model voxelized composites and accurately compute effective material properties. Saputra et al. \cite{saputra2017automatic} proposed an automatic octree meshing algorithm combined with polyhedral SBFEM to efficiently perform stress analysis based on digital images. This integration leverages the hierarchical refinement capabilities of octree structures, enhancing the automatic meshing, adaptivity, and computational efficiency of the SBFEM for complex 3D problems.

In recent studies, the SBFEM has been employed for seepage analysis. Li and Tu \cite{li2012scaled} used the SBFEM to solve steady-state seepage problems with multi-material regions. Liu et al. \cite{liu2018new} presented an iso-geometric SBFEM using nonuniform rational B-splines for the numerical solution of seepage problems in the unbounded domain. Johari and Heydari \cite{johari2018reliability} proposed a stochastic SBFEM for seepage reliability analysis considering spatial variability. Yang et al. \cite{yangNovelSolutionSeepage2022} integrated the SBFEM and polygonal mesh technique to solve steady-state and transient seepage problems. Yan et al. \cite{yan2023psbfem} proposed a PSBFEM approach that integrates quadtree mesh generation from digital images to solve seepage problems. To date, most seepage analyses have been concentrated on two-dimensional (2D) domains. While these 2D approaches have provided valuable insights, they are often limited when it comes to modeling complex, real-world scenarios. Building upon these studies, this work aims to extend the SBFEM approach to 3D seepage problems.

In this study, a polyhedral SBFEM (PSBFEM) framework is proposed for solving 3D seepage problems. The remainder of this paper is organized as follows: Sections \ref{sec:2} and \ref{sec:3} derive the formulation of the SBFEM for 3D seepage problems, including the governing equations and the scaled boundary finite element equations. Section \ref{sec:4} develops a PSBFEM by incorporating Wachspress shape functions and introducing the concept of polyhedronal PSBFEM elements. Section \ref{sec:5} presents the solution procedure for the SBFEM equations, covering stiffness matrix solution, mass matrix solution, transient solution, and free surface seepage solution. Section \ref{sec:6} describes the implementation process for solving seepage problems using the PSBFEM in ABAQUS UEL, explaining how to use the user-defined element to perform relevant calculations. Section \ref{sec:7} provides several numerical examples, such as steady-state, transient seepage and free surface analyses of different structures, to demonstrate the accuracy and efficiency of the proposed method. Finally, Section \ref{sec:8} concludes the study by summarizing the key findings, discussing the limitations of the method, and suggesting directions for future work.

\section{3D seepage equations}
\label{sec:2}
We considered a 3D transient seepage problems in this work, the governing equation without source terms can be expressed as \cite{johari2018reliability} 
\begin{equation}
    \mathbf{L}^\mathrm{T}\mathbf{q}+ S_s\dot{{h}}=0\quad\text{in } \quad \Omega,  \label{eq:gov} 
\end{equation}
where $S_s$ is the specific storage coefficient, $\dot{h}$ is the total head change rate; $\mathbf{q}$ denotes the flux vector. $\mathbf{q}$ can be written as
\begin{equation}
    \mathbf{q}=-\mathbf{k}\mathbf{L}h,
\end{equation}
where $h$ denotes the total hydraulic head, $\mathbf{k}$ is the hydraulic conductivity matrix. The operator $\mathbf{L}$ is the differential operator and can be written as
\begin{equation}
    \mathbf{L}=\begin{Bmatrix}\frac{\partial}{\partial\hat{x}}\\\frac{\partial}{\partial\hat{y}}\\\frac{\partial}{\partial\hat{z}}\end{Bmatrix}.
\end{equation}

By applying the Fourier transform on Eq. (\ref{eq:gov}), the governing equation is transformed into the frequency domain as 
\begin{equation}
\mathbf{L}^\mathrm{T}\mathbf{\tilde{q}}+\mathrm{i}\omega S_s \tilde{h} =0
, \label{eq:threegov}
\end{equation}
where $\mathbf{\tilde{q}}$ and $\tilde{h}$ are the Fourier transform of $\mathbf{q}$ and $h$, respectively. $\omega$ is the frequency. When $\omega=0$, the problems are transformed into the steady-state seepage problem.

\section{Formulation of SBFEM for seepage problems}
\label{sec:3}
\subsection{Scaled boundary transformation of the geometry}
\begin{figure}[]
  \centering
  \includegraphics[width=0.6\textwidth]{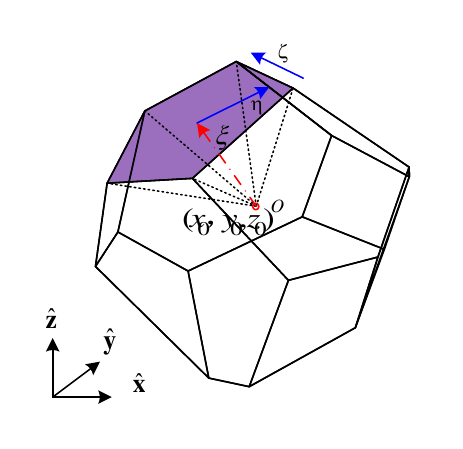}
  \caption{Coordinate formulation $(\xi, \eta, \zeta)$ in a scaled boundary setting for general polyhedral elements.}
  \label{fig:A ployhedral element}
\end{figure}

Fig.~\ref{fig:A ployhedral element} illustrates the three-dimensional scaled boundary coordinate system defined for the polyhedral element. The circumferential coordinates, denoted by $\xi$ and $\eta$, are defined on the surface of the element, while $\zeta$ represents the radial coordinate. The Cartesian coordinates of a point $(\hat{x}, \hat{y}, \hat{z})$ located within the volume sector $\mathrm{V}_{\mathrm{e}}$ can be expressed in terms of the scaled boundary coordinates $(\xi, \eta, \zeta)$ as described in~\cite{song2018scaled}
\begin{subequations}
\begin{align}
    \hat{x}(\xi,\eta,\zeta)=\xi\hat{x}(\eta,\zeta)=\xi\mathbf{N}(\eta,\zeta)\mathbf{\hat{x}}, \\
    \hat{y}(\xi,\eta,\zeta)=\xi\hat{y}(\eta,\zeta)=\xi\mathbf{N}(\eta,\zeta)\mathbf{\hat{y}}, \\
    \hat z(\xi,\eta,\zeta)=\xi\hat z(\eta,\zeta)=\xi\mathbf{N}(\eta,\zeta)\mathbf{\hat{z}},
\end{align}
\end{subequations}
where $\mathbf{\hat{x}}$, $\mathbf{\hat{y}}$ and $\mathbf{\hat{z}}$ are the nodal coordinate vectors of the surface element $\mathrm{S_e}$ in Cartesian coordinates. $\mathbf{N}(\eta,\zeta)$ is the shape function vector
\begin{equation}
    \mathbf{N}(\eta,\zeta)=[ N_1(\eta,\zeta) \quad N_2(\eta,\zeta)  \quad ... \quad N_i(\eta,\zeta)\quad ... \quad N_n(\eta,\zeta)],
\end{equation}
where $N_i(\eta, \zeta)$ are the shape functions and $n$ is the total number of nodes in the polygon. The infinitesimal volume $\mathrm{d}\Omega$ defined in the scaled boundary coordinate system is transformed into the Cartesian coordinate system through the following mapping relation \cite{song2018scaled}
\begin{equation}d\Omega=\xi^{2}|\mathbf{J}|d\xi d\eta d\zeta\end{equation}
where $|\mathbf{J}|$ is the the determinant of the Jacobian matrix $\mathbf{J}(\eta,\zeta)$.

The partial derivatives with respect to the scaled boundary coordinates are related to the partial derivatives in Cartesian coordinates through the following equation:
\begin{equation}\begin{Bmatrix}\frac{\partial}{\partial\xi}\\\frac{\partial}{\partial\eta}\\\frac{\partial}{\partial\zeta}\end{Bmatrix}=\begin{bmatrix}1&0&0\\0&\xi&0\\0&0&\xi\end{bmatrix}\mathbf{J_b}(\eta,\zeta)\begin{Bmatrix}\frac{\partial}{\partial x}\\\frac{\partial}{\partial y}\\\frac{\partial}{\partial z}\end{Bmatrix}, \label{eq:L1}\end{equation}
where $\mathbf{J_b}(\eta,\zeta)$ represents the Jacobian matrix on the boundary $(\xi=1)$, defined as follows:
\begin{equation}\mathbf{J_b}(\eta,\zeta)=\begin{bmatrix}\hat{x}(\eta,\zeta)&\hat{y}(\eta,\zeta)&\hat{z}(\eta,\zeta)\\\hat{x}(\eta,\zeta)_{,\eta}&\hat{y}(\eta,\zeta)_{,\eta}&\hat{z}(\eta,\zeta)_{,\eta}\\\hat{x}(\eta,\zeta)_{,\zeta}&\hat{y}(\eta,\zeta)_{,\zeta}&\hat{z}(\eta,\zeta)_{,\zeta}\end{bmatrix}.\end{equation}

The determinant of $\mathbf{J_b}(\eta,\zeta)$ is expressed as:
\begin{equation}|\mathbf{J_b}|=\hat{x}(\hat{y}_{,\eta}\hat{z}_{,\zeta}-\hat{z}_{,\eta}\hat{y}_{,\zeta})+\hat{y}(\hat{z}_{,\eta}\hat{x}_{,\zeta}-\hat{x}_{,\eta}\hat{z}_{,\zeta})+\hat{z}(\hat{x}_{,\eta}\hat{y}_{,\zeta}-\hat{y}_{,\eta}\hat{x}_{,\zeta}),\end{equation}

For clarity, the argument $(\eta,\zeta)$ has been omitted. From Eq. (\ref{eq:L1}), the inverse relationship is given by:
\begin{equation}\left\{\begin{array}{c}\frac{\partial}{\partial\hat{x}}\\\frac{\partial}{\partial\hat{y}}\\\frac{\partial}{\partial\hat{z}}\end{array}\right\}=\mathbf{J_b}^{-1}\left(\boldsymbol{\eta},\boldsymbol{\zeta}\right)\left\{\begin{array}{c}\frac{\partial}{\partial\xi}\\\frac{1}{\xi}\frac{\partial}{\partial\eta}\\\frac{1}{\xi}\frac{\partial}{\partial\zeta}\end{array}\right\},\label{eq:L2}\end{equation}
where $\mathbf{J_b}^{-1}\left(\eta,\zeta\right)$ denotes the inverse of the Jacobian matrix at the boundary. $\mathbf{J}^{-1}\left(\eta,\zeta\right)$ can be written as
\begin{equation}\left.\mathbf{J}^{-1}\left(\boldsymbol{\eta},\boldsymbol{\zeta}\right)=\frac1{\left|\mathbf{J_b}\right|}\left[\begin{array}{ccc}y_{,\boldsymbol{\eta}}z_{,\boldsymbol{\zeta}}-z_{,\boldsymbol{\eta}}y_{,\boldsymbol{\zeta}}&zy_{,\boldsymbol{\zeta}}-yz_{,\boldsymbol{\zeta}}&yz_{,\boldsymbol{\eta}}-zy_{,\boldsymbol{\eta}}\\\\z_{,\boldsymbol{\eta}}x_{,\boldsymbol{\zeta}}-x_{,\boldsymbol{\eta}}z_{,\boldsymbol{\zeta}}&xz_{,\boldsymbol{\zeta}}-zx_{,\boldsymbol{\zeta}}&zx_{,\boldsymbol{\eta}}-xz_{,\boldsymbol{\eta}}\\\\x_{,\boldsymbol{\eta}}y_{,\boldsymbol{\zeta}}-y_{,\boldsymbol{\eta}}x_{,\boldsymbol{\zeta}}&yx_{,\boldsymbol{\zeta}}-xy_{,\boldsymbol{\zeta}}&xy_{,\boldsymbol{\eta}}-yx_{,\boldsymbol{\eta}}\end{array}\right.\right].\end{equation}


The gradient operator in the Cartesian coordinate system can be transformed into the scaled boundary coordinate system ($\xi, \eta, \zeta$) as follows:
\begin{equation}\mathbf{L}=\mathbf{b}_{1}\left(\eta,\zeta\right)\frac{\partial}{\partial\xi}+\frac{1}{\xi}\left(\mathbf{b}_{2}\left(\eta,\zeta\right)\frac{\partial}{\partial \eta}+\mathbf{b}_{3}\left(\eta,\zeta\right)\frac{\partial}{\partial\zeta}\right),\label{eq:L3}\end{equation}
where
\begin{equation}
    \mathbf{b}_1\left(\eta,\zeta\right)=\left.\frac1{|\mathbf{J_b}|}\left[\begin{array}{c}y_{,\eta} z_{,\zeta}-z_{,\eta} y_{,\zeta}\\\\z_{,\eta} x_{,\zeta}-x_{,\eta} z_{,\zeta}\\\\x_{,\eta} y_{,\zeta}-y_{,\eta} x_{,\zeta}\end{array}\right.\right],
\end{equation}
\begin{equation}\mathbf{b}_2\left(\eta,\zeta\right)=\left.\frac{1}{|\mathbf{J_b}|}\left[\begin{array}{c}zy_{,\zeta}-yz_{,\zeta}\\\\xz_{,\zeta}-zx_{,\zeta}\\\\yx_{,\zeta}-xy_{,\zeta}\end{array}\right.\right],\end{equation}
\begin{equation}\mathbf{b}_3\left(\eta,\zeta\right)=\left.\frac{1}{|\mathbf{J_b}|}\left[\begin{array}{c}yz_{,\eta}-zy_{,\eta}\\\\zx_{,\eta}-xz_{,\eta}\\\\xy_{,\eta}-yx_{,\eta}\end{array}\right.\right].\end{equation}

\subsection{Seepage field}
The hydraulic head of any point, $\tilde{h}(\xi,\eta,\zeta)$, in the SBFEM coordinates can be expressed as
\begin{equation}
    \tilde{h}(\xi,\eta,\zeta)=\mathbf{N}(\eta,\zeta)\tilde{h}(\xi), \label{eq:L4}
\end{equation}
where $\tilde{h}(\xi)$ is the nodal hydraulic head vector and $\mathbf{N}(\eta,\zeta)$ represents the matrix of shape functions.

By using Eqs. (\ref{eq:L4}) and (\ref{eq:L3}), the partial derivatives of the hydraulic head field are obtained as follows:
\begin{equation}\mathbf{L}=\mathbf{b}_{1}\left(\eta,\zeta\right)\frac{\partial \tilde{h}}{\partial\xi}+\frac{1}{\xi}\left(\mathbf{b}_{2}\left(\eta,\zeta\right)\frac{\partial \tilde{h}}{\partial \eta}+\mathbf{b}_{3}\left(\eta,\zeta\right)\frac{\partial \tilde{h}}{\partial\zeta}\right), \label{eq:L5}\end{equation}

For simplicity, Eq. (\ref{eq:L5}) can be rewritten by substituting the hydraulic head expression from Eq. (\ref{eq:L4}), as follows
\begin{equation}
\mathbf{L}=\mathbf{B}_1\left(\eta,\zeta\right)\tilde{h}_{,\xi}+\frac1\xi \mathbf{B}_2\left(\eta,\zeta\right) \tilde{h},
\end{equation}
where $\mathbf{B}_1\left(\eta,\zeta\right)$ and $\mathbf{B}_2\left(\eta,\zeta\right)$ can be written as
\begin{equation}
\mathbf{B}_1\left(\eta,\zeta\right)=\mathbf{b}_{1}\left(\eta,\zeta\right)\mathbf{N}(\eta,\zeta),
\end{equation}
\begin{equation}
\mathbf{B}_2\left(\eta,\zeta\right)=\mathbf{b}_{2}\left(\eta,\zeta\right)\mathbf{N}(\eta,\zeta)_{,\eta}+\mathbf{b}_{3}\left(\eta,\zeta\right)\mathbf{N}(\eta,\zeta)_{,\zeta}.
\end{equation}

The flux, $\tilde{\mathbf{q}}(\xi,\eta,\zeta)$, can be expressed in the SBFEM coordinate system ($\xi, \eta, \zeta$) as follows:
\begin{equation}\tilde{\mathbf{q}}(\xi,\eta)=-\mathbf{k}\Big(\mathbf{B_1}(\eta,\zeta) \tilde{h}(\xi)_{,\xi}+\frac{1}{\xi}\mathbf{B_2}(\eta,\zeta) \tilde{h}(\xi) \Big),\end{equation}

\subsection{Scaled boundary finite element equation}
By applying the method of weighted residuals \cite{songScaled1999}, Eq. (\ref{eq:threegov}) can be transformed into
\begin{equation}
\begin{aligned}&\int_\Omega w\mathbf{b_1}^\mathrm{T}\mathbf{\tilde{q},_\xi}\mathrm{~d}\Omega+\int_\Omega w\frac{1}{\xi}(\mathbf{b_2}^\mathrm{T}\mathbf{\tilde{q},_\eta}+\mathbf{b_3}^\mathrm{T}\mathbf{\tilde{q},_\zeta})\mathrm{~d}\Omega\\&
+\mathrm{i}\omega\int_\Omega w\rho c\tilde{h}\mathrm{d}\Omega=0,\end{aligned}\label{eq:weighted gov}
\end{equation}
where $w = w(\xi, \eta, \zeta)$ is the weighting function. Following the approach outlined by \cite{songScaled1999}, Eq. (\ref{eq:weighted gov}) is further simplified as follows:
\begin{align}
    &\mathbf{E}_0 \xi^2 \tilde{h}(\xi)_{,\xi\xi} 
    + \left( 2\mathbf{E}_0 - \mathbf{E}_1 + \mathbf{E}_1^{\mathrm{T}} \right) \xi \tilde{h}(\xi)_{,\xi} \notag \\
    &+ \left( \mathbf{E}_1^\mathrm{T} - \mathbf{E}_2 \right) \tilde{h}(\xi) 
    - \mathbf{M}_0 \mathrm{i}\omega \xi^2 \tilde{h}(\xi) = \xi \mathbf{F}(\xi),  
    \label{eq:mainequation1}
\end{align}

The coefficient matrices $\mathbf{E_0}$, $\mathbf{E_1}$, $\mathbf{E_2}$, and $\mathbf{M_0}$ of the entire element are assembled from the corresponding matrices $\mathbf{E_{0}^{e}}$, $\mathbf{E_{1}^{e}}$, $\mathbf{E_{2}^{e}}$, and $\mathbf{M_{0}^{e}}$ of each surface element. The coefficient matrices $\mathbf{E_{0}^{e}}$, $\mathbf{E_{1}^{e}}$, $\mathbf{E_{2}^{e}}$, and $\mathbf{M_{0}^{e}}$ for a surface element $\mathrm{S_e}$ can be written as
\begin{equation}\mathbf{E}_0^\mathrm{e}=\int_{\mathbf{S}_\mathrm{e}}\mathbf{B}_1^\mathrm{T}\mathbf{kB}_1\left|\mathbf{J}_\mathrm{b}\right|\mathrm{d}\eta\mathrm{d}\zeta,\label{eq:E0}\end{equation}
\begin{equation}\mathbf{E}_1^\mathrm{e}=\int_{\mathbf{S}_\mathrm{e}}\mathbf{B}_2^\mathrm{T}\mathbf{kB}_1\left|\mathbf{J}_\mathrm{b}\right|\mathrm{d}\eta\mathrm{d}\zeta,\end{equation}
\begin{equation}\mathbf{E}_2^\mathrm{e}=\int_{\mathrm{S}_\mathrm{e}}\mathbf{B}_2^\mathrm{T}\mathbf{kB}_2\left|\mathbf{J}_\mathrm{b}\right|\mathrm{d}\eta\mathrm{d}\zeta,\end{equation}
\begin{equation}\mathbf{M_0^e}=\int_{\mathrm{S_e}}\mathbf{N^T}S_s\mathbf{N}\left|\mathbf{J_b}\right|\mathrm{d}\eta\mathrm{d}\zeta,\label{eq:M0}\end{equation}
where $\mathbf{k}$ is the matrix of permeability coefficients.
\section{Polyhedronal SBFEM element}
\label{sec:4}
\subsection{Polyhedral element construction}
\begin{figure}[H]
  \centering
  \includegraphics[width=1.0\textwidth]{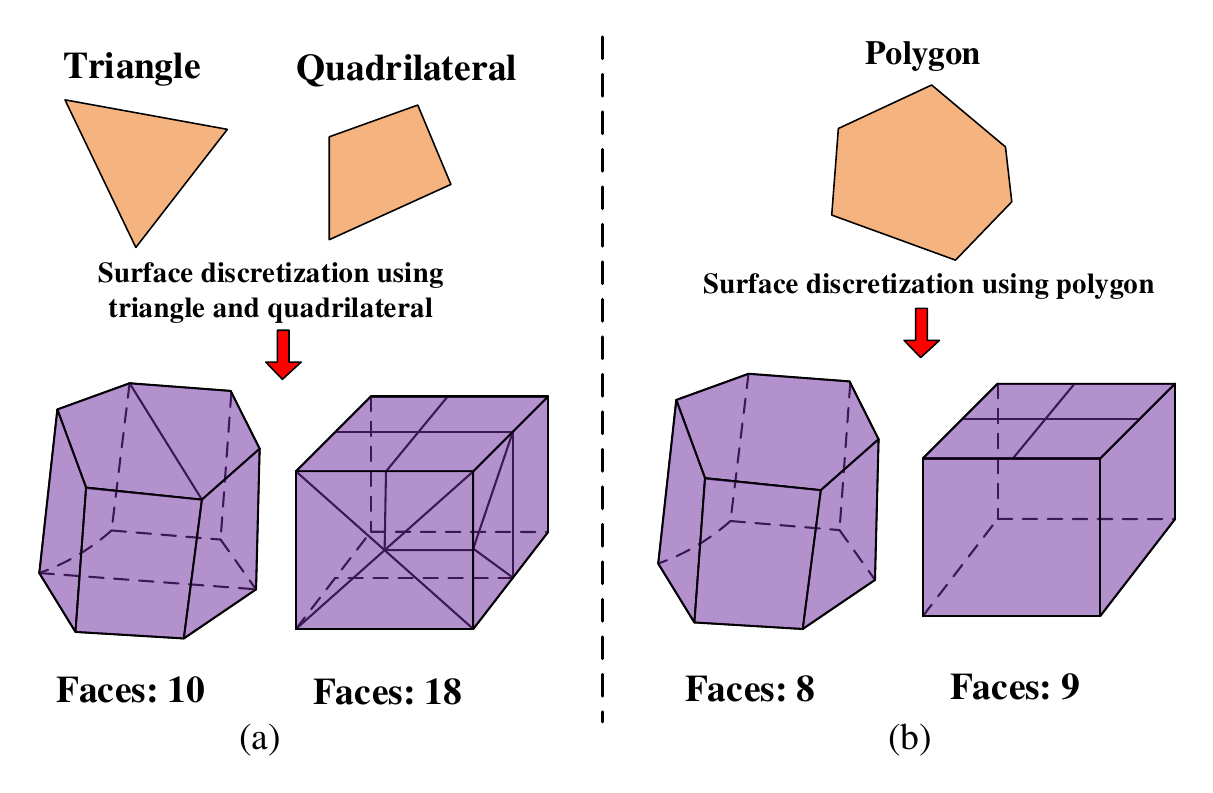}
  \caption{Polyhedral elements construction based on polygonal surfaces. (a) conventional surface discretization method; (b) surface discretization using polygon.}
  \label{fig:polyhedral_element}
\end{figure}

As illustrated in Fig. \ref{fig:polyhedral_element}(a), boundary surfaces are discretized using triangular and quadrilateral elements \cite{yangDevelopment2020,ya2021open}. As a result, the polyhedra processed by traditional 3D SBFEM exhibit relatively complex topologies. To address this complexity, this work introduces polygonal discretization techniques that simplify the topological structure of the polyhedra, thereby reducing the number of element faces, as shown in Fig \ref{fig:polyhedral_element}(b). From the figures, it is clear that the polyhedra constructed using polygonal discretization significantly reduce the number of element faces. This reduction not only enhances mesh efficiency but also provides clearer visualization.

\begin{figure}[H]
  \centering
  \includegraphics[width=1.0\textwidth]{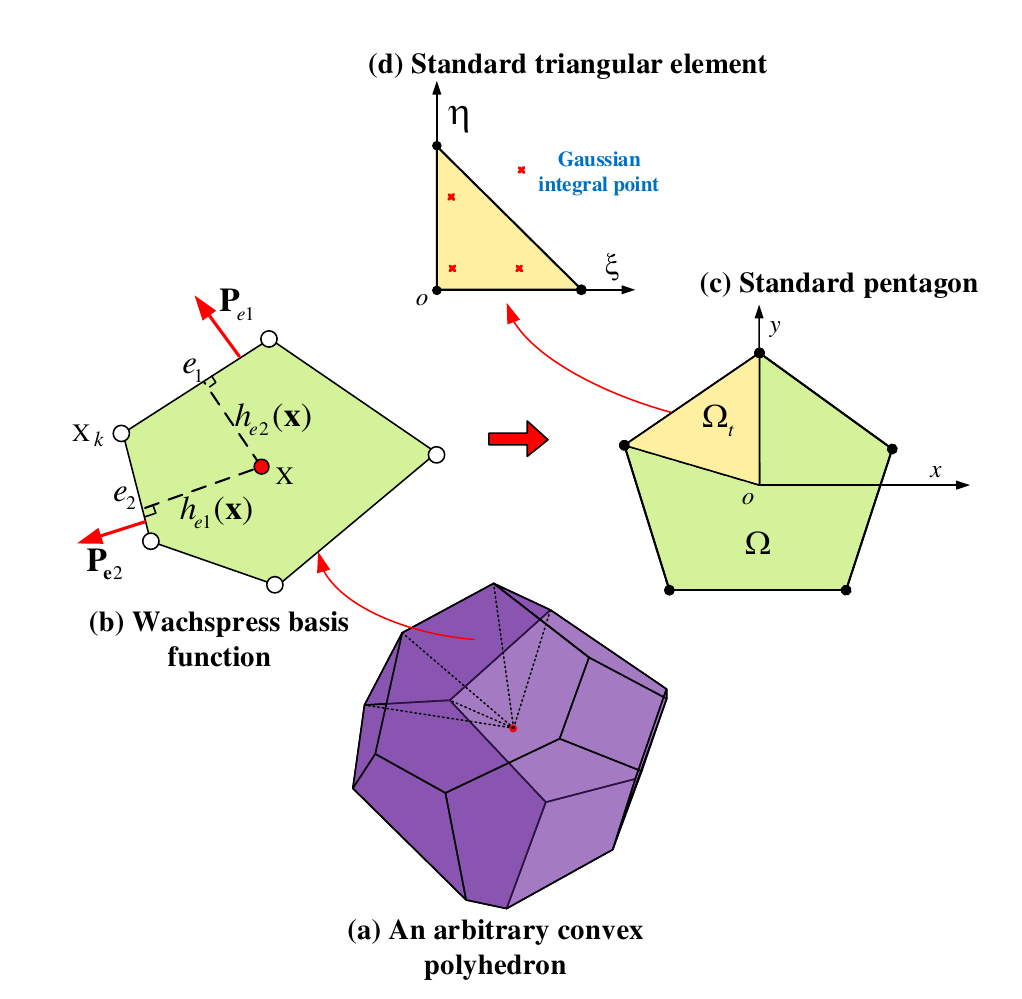}
  \caption{Numerical integration techniques for an arbitrary convex polyhedron; (a) an arbitrary convex polyhedron; (b) wachspress basis function; (c) standard triangular element; (d) standard pentagon.}
  \label{fig:Numerical integration technique}
\end{figure}

In this work, the Wachspress shape functions are introduced to construct polygonal elements. Wachspress~\cite{wachspress2006rational} proposed rational basis functions for convex polygonal elements based on principles from projective geometry. These functions preserve nodal interpolation and ensure linearity along the edges by employing algebraic representations of the polygon boundaries. Warren~\cite{warren2003uniqueness, warren2007barycentric} later extended this formulation to convex polyhedra, as illustrated in Fig.~\ref{fig:Numerical integration technique}(b).

Let \( P \subset \mathbb{R}^2 \) be a simple convex polygon with a set of edges \( E \) and vertices \( V \). For each edge \( e \in E \), let \( \mathbf{n}_e \) denote its unit outward normal, and \( h_e(\mathbf{x}) \) the perpendicular distance from a point \( \mathbf{x} \in P \) to edge \( e \), computed as:
\begin{equation}
h_e(\mathbf{x}) = (\mathbf{v}_e - \mathbf{x}) \cdot \mathbf{n}_e,
\end{equation}
where \( \mathbf{v}_e \) is any point on edge \( e \). Define the scaled normal vector as:
\begin{equation}
\mathbf{p}_e(\mathbf{x}) = \frac{\mathbf{n}_e}{h_e(\mathbf{x})}.
\end{equation}

For a vertex \( \mathbf{v} \in V \) adjacent to two edges \( e_1, e_2 \), the unnormalized weight function is given by:
\begin{equation}
w_{\mathbf{v}}(\mathbf{x}) = \frac{
\det(\mathbf{p}_{e_1}(\mathbf{x}), \mathbf{p}_{e_2}(\mathbf{x}))
}{
h_{e_1}(\mathbf{x})\, h_{e_2}(\mathbf{x})
}.
\end{equation}

Finally, the Wachspress shape function at vertex \( \mathbf{v} \) is:
\begin{equation}
N(\mathbf{x}) = \frac{w_{\mathbf{v}}(\mathbf{x})}{\sum_{\mathbf{u} \in V} w_{\mathbf{u}}(\mathbf{x})}.
\end{equation}

\subsection{Numerical integration technique for polygonal surface element $S_e$}
To compute coefficient matrices $\mathbf{E_{0}^{e}}$, $\mathbf{E_{1}^{e}}$, $\mathbf{E_{2}^{e}}$, and $\mathbf{M_{0}^{e}}$ for a polygonal surface $\mathrm{S_e}$, a triangulation-based numerical integration approach is employed in this study. The surface element is initially divided into a set of non-overlapping triangular subdomains. Gaussian quadrature is then applied to each triangle to perform the integration, as presented in Fig. \ref{fig:Numerical integration technique}(d). The results from all the subdomains are subsequently combined to construct coefficient matrices of the entire polygonal surface.

Consider a polygonal surface $\Omega$ with $m$ edges, as shown in Fig. \ref{fig:Numerical integration technique}(c). This domain is subdivided into $m$ triangular subdomains $\Omega_t$, which can be constructed by either introducing an internal node or by connecting all edges to the centroid $\mathbf{x}_c$, such that
\begin{equation}
\Omega = \bigcup_{t=1}^{m} \Omega_t,
\end{equation}
where $\Omega_t$ represents the $t$-th triangular subdomain. For each triangular subdomain $\Omega_t$, the numerical integration is performed using Gaussian quadrature, expressed as
\begin{equation}
\int_{\Omega} f(\mathbf{x}) \, d\mathbf{x} = \sum_{t=1}^{m} \sum_{i=1}^{n_g} w_i f(\mathbf{x}_i),
\end{equation}
where $\mathbf{x}_i$ are the Gaussian integration points within the triangle $\Omega_t$, $w_i$ are the corresponding integration weights, and $n_g$ is the total number of quadrature points in each subdomain.

Given the shape functions $N_k(\mathbf{x})$ of the polygonal surface, using Gaussian quadrature over triangular subdomains, the final expression of coefficient matrices becomes:
\begin{equation}
\mathbf{E}_0^\mathrm{e}=\sum_{t=1}^{m} \sum_{i=1}^{n_g} w_i\mathbf{B}_1(\mathbf{x}_i)^\mathrm{T}\mathbf{kB}_1(\mathbf{x}_i)\left|\mathbf{J}_\mathrm{b}\right|,
\end{equation}
\begin{equation}
\mathbf{E}_1^\mathrm{e}=\sum_{t=1}^{m} \sum_{i=1}^{n_g} w_i\mathbf{B}_2(\mathbf{x}_i)^\mathrm{T}\mathbf{kB}_1(\mathbf{x}_i)\left|\mathbf{J}_\mathrm{b}\right|,
\end{equation}
\begin{equation}
\mathbf{E}_2^\mathrm{e}=\sum_{t=1}^{m} \sum_{i=1}^{n_g} w_i\mathbf{B}_2(\mathbf{x}_i)^\mathrm{T}\mathbf{kB}_2(\mathbf{x}_i)\left|\mathbf{J}_\mathrm{b}\right|,
\end{equation}
\begin{equation}
\mathbf{M_0^e}=\sum_{t=1}^{m} \sum_{i=1}^{n_g} w_i\mathbf{N^T}(\mathbf{x}_i)S_s\mathbf{N}(\mathbf{x}_i)\left|\mathbf{J_b}\right|.
\end{equation}

\section{Solution procedure}
\label{sec:5}
\subsection{Stiffness matrix solution}
In this work, the term $\mathbf{F}(\boldsymbol{\xi})$ on the right-hand side of Eq.~(\ref{eq:mainequation1}) denotes the contribution of the normal fluid flux across the lateral boundary of the domain. When the lateral boundary is impermeable or the domain is closed, no normal flow occurs through the side faces, leading to $\mathbf{F}(\boldsymbol{\xi}) = 0$ \cite{yu2021scaled}. Under this impermeable boundary condition assumption, the term is neglected in the current analysis. Consequently, the steady-state hydraulic head distribution within the scaled boundary finite element framework can be obtained by substituting the eigenvalue $\omega$ into the governing equation, resulting in the following expression
\begin{equation}\mathbf{E}_0\xi^2\tilde{h}(\xi)_{,\xi\xi}+\left(2\mathbf{E}_0-\mathbf{E}_1+\mathbf{E}_1^\mathrm {h}\right)\xi\tilde{h}(\xi)_{,\xi}+\left(\mathbf{E}_1^h-\mathbf{E}_2\right)\tilde{h}(\xi)=\mathbf{0}.\end{equation}

By introducing the variable, 
\begin{equation}\mathbf{X}(\xi)=\begin{bmatrix}\xi^{0.5}\tilde{\mathbf{h}}(\xi)\\\xi^{-0.5}\tilde{\mathbf{Q}}(\xi)\end{bmatrix},\end{equation}

The SBFEM equation is then reformulated as: 
\begin{equation}\xi\mathbf{X}(\xi)_{,\xi}=\mathbf{Z}_\mathrm{p}\mathbf{X}(\xi),\end{equation}
where the coefficient matrix $\mathbf{Z}_\mathrm{p}$ is a Hamiltonian matrix, defined as: 
\begin{equation}\mathbf{Z_p}=\begin{bmatrix}-\mathbf{E}_0^{-1}\mathbf{E}_1^\mathrm{T}+0.5\mathbf{I}&\mathbf{E}_0^{-1}\\\mathbf{E}_2-\mathbf{E}_1\mathbf{E}_0^{-1}\mathbf{E}_1^\mathrm{T}&\mathbf{E}_1\mathbf{E}_0^{-1}-0.5\mathbf{I}\end{bmatrix}.\label{eq:zp}\end{equation}

The eigenvalue decomposition of $\mathbf{Z}_\mathrm{p}$ is expressed as 
\begin{equation}\mathbf{Z}_\mathrm{p}\begin{bmatrix}\mathbf{\Phi}_\mathrm{h1}&&\mathbf{\Phi}_\mathrm{h2}\\\\\mathbf{\Phi}_\mathrm{q1}&&\mathbf{\Phi}_\mathrm{q2}\end{bmatrix}=\begin{bmatrix}\mathbf{\Phi}_\mathrm{h1}&&\mathbf{\Phi}_\mathrm{h2}\\\\\mathbf{\Phi}_\mathrm{q1}&&\mathbf{\Phi}_\mathrm{q2}\end{bmatrix}\begin{bmatrix}\boldsymbol{\Lambda}^+&&0\\\\0&&\boldsymbol{\Lambda}^-\end{bmatrix},\label{eq:eigen decomp}\end{equation}

Here, $\mathbf{\Phi}_{\mathrm{h1}}$ and $\mathbf{\Phi}_{\mathrm{q1}}$ represent the eigenvectors associated with the diagonal matrix $\boldsymbol{\Lambda}^+$, while $\mathbf{\Phi}_{\mathrm{h2}}$ and $\mathbf{\Phi}_{\mathrm{q2}}$ correspond to those of $\boldsymbol{\Lambda}^-$. The matrices $\boldsymbol{\Lambda}^+$ and $\boldsymbol{\Lambda}^-$ consist of eigenvalues with positive and negative real parts, respectively. Collectively, $\mathbf{\Phi}_{\mathrm{h}}$ and $\mathbf{\Phi}_{\mathrm{q}}$ characterize the modal hydraulic head and flux components. The solution vector $\mathbf{X}(\xi)$ can thus be expressed as

\begin{equation}\mathbf{X}(\xi)=\begin{bmatrix}\mathbf{\Phi}_\mathrm{h1}&&\mathbf{\Phi}_\mathrm{h2}\\\\\mathbf{\Phi}_{\mathrm{q1}}&&\mathbf{\Phi}_\mathrm{q2}\end{bmatrix}\begin{bmatrix}\xi^{\boldsymbol{\Lambda}^+}&&\mathbf{0}\\\\\mathbf{0}&&\xi^{\boldsymbol{\Lambda}^-}\end{bmatrix}\begin{Bmatrix}\mathbf{c}^{\mathrm{\mathrm{n}}}\\\\\mathbf{c}^{\mathrm{p}}\end{Bmatrix},\end{equation}
where the vectors $\mathbf{c}^{\mathrm{n}}$ and $\mathbf{c}^{\mathrm{p}}$ denote the integration constants associated with the eigenvalue matrices $\boldsymbol{\Lambda}^+$ and $\boldsymbol{\Lambda}^-$, respectively. These constants are determined by enforcing the relevant boundary conditions. Accordingly, the solutions for the hydraulic head $\mathbf{\tilde{h}}(\xi)$ and flux $\mathbf{\tilde{q}}(\xi)$ can be written as
\begin{equation}\mathbf{\tilde{h}}(\xi)=\mathbf{\Phi}_{\mathrm{h1}}\xi^{-\boldsymbol{\Lambda}^+-0.5\mathbf{I}}\mathbf{c}^{\mathrm{n}},\end{equation}
\begin{equation}\mathbf{\tilde{q}(\xi)=\Phi}_{\mathrm{q}1}\xi^{-\boldsymbol{\Lambda}^++0.5\mathbf{I}}\mathbf{c}^{\mathrm{n}}.\end{equation}
where the integration constants $\mathbf{c}^{\mathrm{n}}$ is determined from the boundary conditions. 

The relation between the nodal hydraulic head functions and nodal flux functions is derived by eliminating the integration constants, as expressed below:
\begin{equation}\mathbf{\tilde{q}}(\xi)=\mathbf{\Phi}_{\mathrm{q1}}\mathbf{\Phi}_{\mathrm{h1}}^{-1}\xi\mathbf{\tilde{h}}(\xi).\end{equation}

At the boundary $(\xi=1)$, the nodal flux vector is denoted by $\mathbf{Q} = \mathbf{\tilde{q}}(\xi=1)$, and the nodal hydraulic head vector is given by $\mathbf{h} = \mathbf{\tilde{h}}(\xi=1)$. The relationship between the flux and hydraulic head vectors can be written as
\[
\mathbf{Q} = \mathbf{K} \mathbf{h}.
\]
Accordingly, the stiffness matrix for the subdomain is formulated as:
\begin{equation}\mathbf{K}=\mathbf{\Phi}_{\mathrm{q1}}\mathbf{\Phi}_{\mathrm{h1}}^{-1}.\label{eq:K}\end{equation}

\subsection{Mass matrix solution}
The mass matrix of a volume element is formulated as \cite{yan2025polyhedral,song2018scaled}
\begin{equation}
\mathbf{M} = \mathbf{\Phi}_{\mathrm{h}1}^{-\mathrm{T}} \left( \int_0^1 \xi^{\boldsymbol{\Lambda}^+} \mathbf{m}_0 \, \xi^{\boldsymbol{\Lambda}^+} \, \xi \, \mathrm{d}\xi \right) \mathbf{\Phi}_{\mathrm{h}1}^{-1}, \label{eq:Massequation}
\end{equation}
where the coefficient matrix $\mathbf{m}_0$ is defined by
\begin{equation}
\mathbf{m}_0 = \mathbf{\Phi}_{\mathrm{h}1}^\mathrm{T} \mathbf{M}_0 \mathbf{\Phi}_{\mathrm{h}1}.
\end{equation}

By expressing the integral in Eq.~(\ref{eq:Massequation}) in matrix form, the mass matrix can be equivalently written as
\begin{equation}
\mathbf{M} = \mathbf{\Phi}_{\mathrm{h}1}^{-\mathrm{T}} \mathbf{m} \mathbf{\Phi}_{\mathrm{h}1}^{-1}, \label{eq:M}
\end{equation}
where
\begin{equation}
\mathbf{m} = \int_0^1 \xi^{\boldsymbol{\Lambda}^+} \mathbf{m}_0 \, \xi^{\boldsymbol{\Lambda}^+} \, \xi \, \mathrm{d}\xi.
\end{equation}

Each entry of the matrix $\mathbf{m}$ can be evaluated analytically as
\begin{equation}
m_{ij} = \frac{m_{0ij}}{\lambda_{ii}^+ + \lambda_{jj}^+ + 2},
\end{equation}
where $m_{0ij}$ is the $(i,j)$-th entry of $\mathbf{m}_0$, and $\lambda_{ii}^+$, $\lambda_{jj}^+$ are the diagonal entries of the eigenvalue matrix $\boldsymbol{\Lambda}^+$.

\subsection{Transient solution}
The governing equation describing the nodal hydraulic head within a bounded domain can be formulated in the time domain as:
\begin{equation}
\mathbf{K}\mathbf{h}(t) + \mathbf{M}\dot{\mathbf{h}}(t) = \mathbf{Q}(t),
\label{eq:L6}
\end{equation}
where $\mathbf{h}(t)$ denotes the time-dependent nodal hydraulic head, assumed to be continuously differentiable. Due to the complexity of directly solving Eq.~(\ref{eq:L6}) analytically in the time domain, a numerical approach is adopted. In this study, the backward difference scheme~\cite{zienkiewicz1989finite} is employed to approximate the temporal derivative. The simulation time is discretized into uniform intervals, and the solution is obtained sequentially using the prescribed initial conditions. Intermediate values of the hydraulic head can subsequently be determined via interpolation.

Within a time increment $\left[t, t + \Delta t\right]$, the temporal derivative of the hydraulic head is approximated as:
\begin{equation}
\dot{\mathbf{h}}(t) \approx \frac{\Delta \mathbf{h}}{\Delta t} = \frac{\mathbf{h}^{t+\Delta t} - \mathbf{h}^t}{\Delta t}.
\label{eq:L7}
\end{equation}

Substituting Eq.~(\ref{eq:L7}) into Eq.~(\ref{eq:L6}) results in the following discretized system at time $t + \Delta t$:
\begin{equation}
\left(\mathbf{K}^{t+\Delta t} + \frac{\mathbf{M}^{t+\Delta t}}{\Delta t} \right)\mathbf{h}^{t+\Delta t} = \mathbf{Q}^{t+\Delta t} + \frac{\mathbf{M}^{t+\Delta t}}{\Delta t} \mathbf{h}^t.
\end{equation}

\section{Implementation}
\label{sec:6}
\subsection{Overview framework}
\begin{figure}[H]
  \centering
  \includegraphics[width=0.9\textwidth]{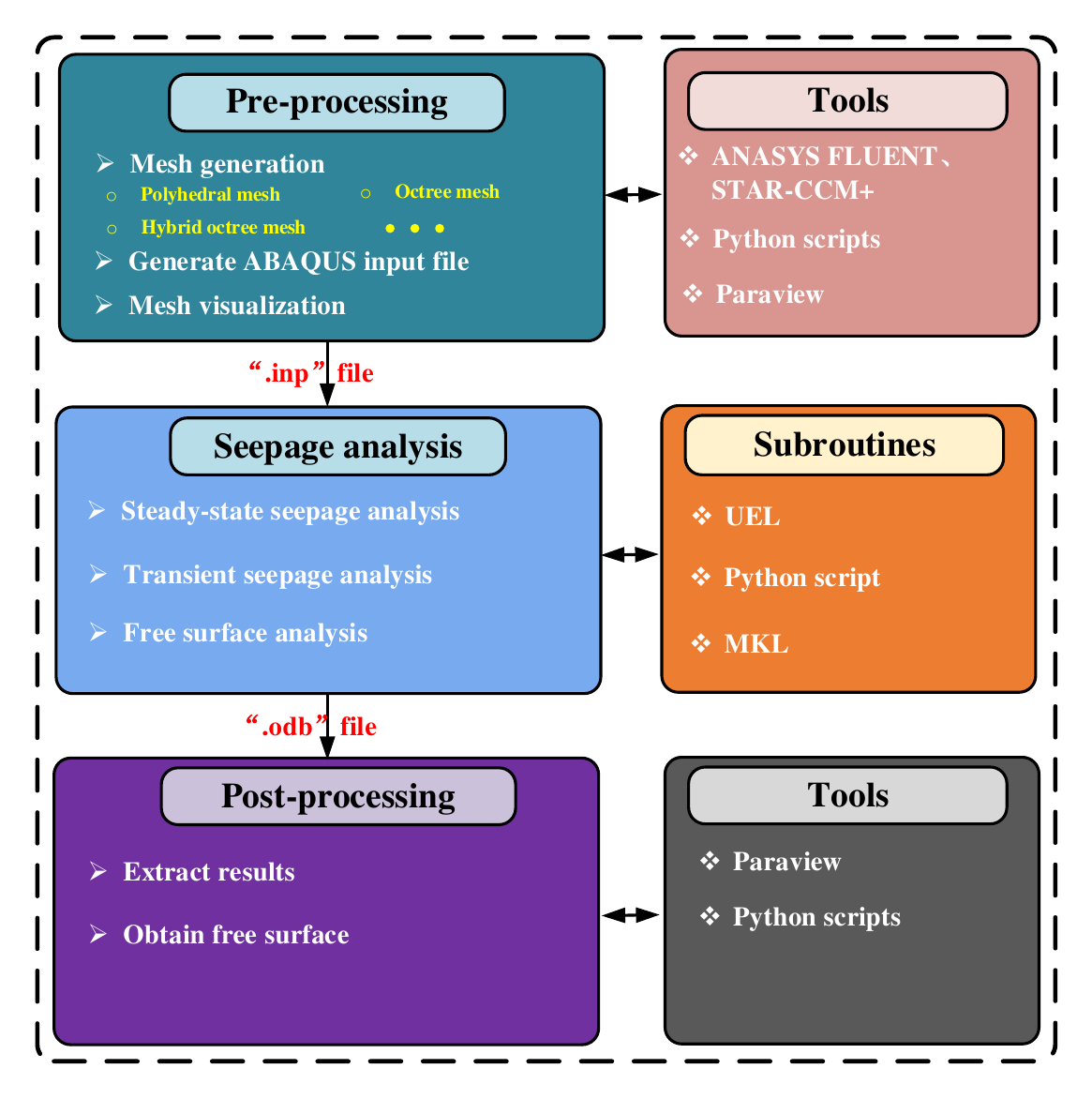}
  \caption{The framework of PSBFEM for 3D seepage analysis.}
  \label{fig:framework}
\end{figure}

A three-dimensional seepage analysis framework is established based on the PSBFEM, as illustrated in Fig.~\ref{fig:framework}. The framework is structured into three main components: pre-processing, numerical analysis, and post-processing. In the pre-processing phase, polyhedral and octree meshes are first generated using commercial CFD software such as ANSYS Fluent~\cite{ANSYSFluent} and STAR-CCM+~\cite{STARCCM}. A customized Python script is subsequently employed to convert the mesh information into the ABAQUS input format (.inp).  

The numerical analysis module incorporates steady-state, transient, and free surface seepage simulations. Both steady and transient seepage problems are solved through a user-defined element (UEL) subroutine written in FORTRAN within the ABAQUS environment. For free surface simulations, a Python-based algorithm is coupled with the UEL framework to capture the evolving seepage boundary. Implementation details are discussed in Sections~\ref{subsec:UEL} and~\ref{subsec:free surface}.  

Given the limitations of ABAQUS/CAE in visualizing UEL elements, the simulation results are extracted from the output database (.odb) and post-processed using ParaView~\cite{ParaViewGuide} for visualization and interpretation.

\subsection{Implementation of the UEL}
\label{subsec:UEL}
The procedure for conducting both steady-state and transient seepage analyses is outlined in Algorithm \ref{alg:1}, and has been implemented in ABAQUS through the User Element (UEL) subroutine interface. Within this framework, the UEL is tasked with computing and assembling the element-level contributions to the global residual vector (RHS) and stiffness matrix (AMATRX). In the context of steady-state seepage analysis, the corresponding expressions for AMATRX and RHS are formulated as follows:
\begin{equation}
\mathrm{AMATRX}=\mathbf{K},\label{eq:amatrx-steady-state}
\end{equation}
\begin{equation}
    \mathrm{RHS}=-\mathbf{K}\mathbf{U},\label{eq:rhs-steady-state}
\end{equation}
where $\mathbf{U}$ is the nodal hydraulic head vector.

For the transient seepage, AMATRX and RHS are defined as follows:
\begin{equation}
\mathrm{AMATRX}=\mathbf{K}^{t+\Delta t}+\frac{\mathbf{M}^{t+\Delta t}}{\Delta t}, \label{eq:amatrx-tranisent-state}
\end{equation}
\begin{equation}
    \mathrm{RHS}=-\mathbf{K}^{t+\Delta t}\mathbf{U}^{t+\Delta t}-\frac{\mathbf{M}^{t+\Delta t}}{\Delta t}\Big(\mathbf{U}^{t+\Delta t}-\mathbf{U}^{t}\Big). \label{eq:rhs-tranisent-state}
\end{equation}

\begin{figure}[H]
  \centering
  \includegraphics[width=0.8\textwidth]{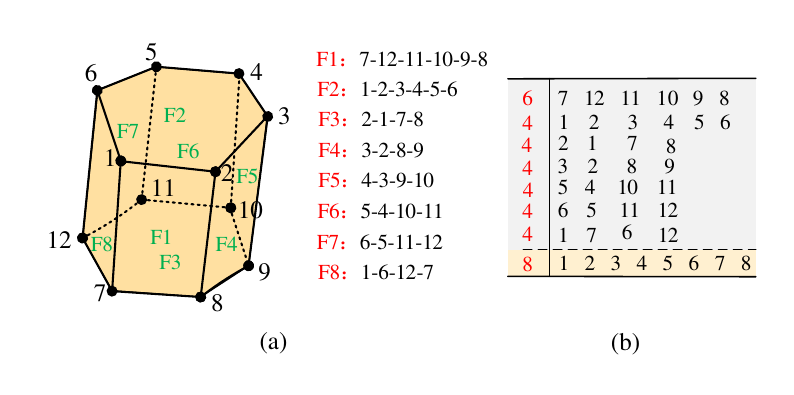}
  \caption{An example of face numbering for polyhedral elements; (a) octahedral element; (b) data structure of element.}
  \label{fig:Element_face}
\end{figure}

To accommodate arbitrarily shaped elements, the formulation adopts a face-by-face construction strategy, offering flexibility for complex polyhedral geometries. To ensure a positive Jacobian determinant, the normal vector of each face must point outward from the element domain. As shown in Fig.~\ref{fig:Element_face}(a), an example octahedral element is presented with labeled face numbers and corresponding node connectivity. Fig.~\ref{fig:Element_face}(b) illustrates the data structure of a polyhedral element: the gray section stores node connectivity for each face, with red labels indicating the number of nodes per face; the orange section stores face connectivity for each element, with red labels denoting the number of faces per element.

\begin{algorithm}
    \caption{Solving 3D Seepage Problems Using PSBFEM}
    \label{alg:1}
    \textbf{Input:} Node and element data, material properties, initial nodal hydraulic head $\mathbf{u}_t$\\
    \textbf{Output:} Updated nodal hydraulic head $\mathbf{u}_{t+1}$
    \begin{algorithmic}[1]
        \While{ABAQUS not converged}
            \State Solve intermediate hydraulic head $\mathbf{u}_{t+1}^{k}$
            \For{each element in AllEle} \Comment{Loop over all elements}
                \State Construct element geometry
                \State Compute element centroid
                \State Compute coefficient matrices $\mathbf{E}_0^e$, $\mathbf{E}_1^e$, $\mathbf{E}_2^e$, and $\mathbf{M}_0^e$ using Eqs.~(\ref{eq:E0})--(\ref{eq:M0})
                \State Assemble matrix $\mathbf{Z}_\mathrm{p}$ using Eq.~(\ref{eq:zp})
                \State Solve eigenvalue decomposition of $\mathbf{Z}_\mathrm{p}$ using Eq.~(\ref{eq:eigen decomp})
                \State Compute stiffness matrix $\mathbf{K}$ and mass matrix $\mathbf{M}$ via Eqs.~(\ref{eq:K}) and (\ref{eq:M})
                \If{$\texttt{lflags}(1) = 62$ or $\texttt{lflags}(1) = 63$}
                    \State Update stiffness matrix $\mathbf{AMATRX}$ and residual vector $\mathbf{RHS}$ using Eqs.~(\ref{eq:amatrx-steady-state}) and (\ref{eq:rhs-steady-state})
                \EndIf
                \If{$\texttt{lflags}(1) = 64$ or $\texttt{lflags}(1) = 65$}
                    \State Update stiffness matrix $\mathbf{AMATRX}$ and residual vector $\mathbf{RHS}$ using Eqs.~(\ref{eq:amatrx-tranisent-state}) and (\ref{eq:rhs-tranisent-state})
                \EndIf
            \EndFor
            \State Update iteration index: $k \gets k + 1$
            \State Set $\mathbf{u}_{t+1} \gets \mathbf{u}_{t+1}^{k}$
        \EndWhile
    \end{algorithmic}
\end{algorithm}

\subsection{Free surface seepage solution}
\label{subsec:free surface}
\begin{figure}[H]
  \centering
  \includegraphics[width=0.8\textwidth]{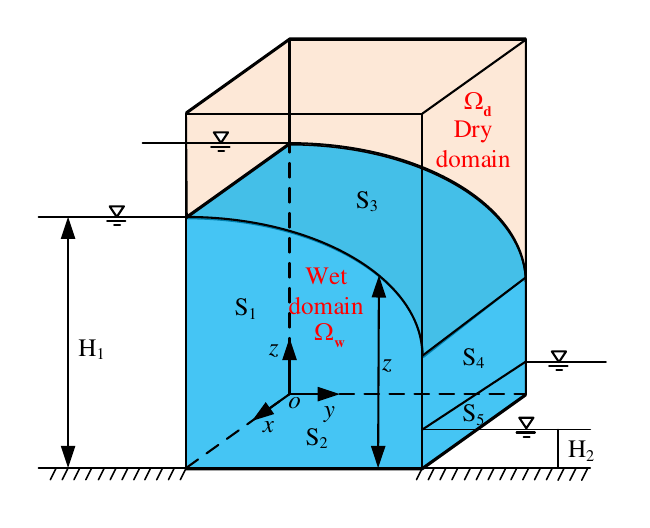}
  \caption{Geometry and boundary conditions of a soil dam.}
  \label{fig:Free surfaces}
\end{figure}

For the seepage problem of free surface in the cubic domain shown in Fig. \ref{fig:Free surfaces}, the problem is divided by the free surface into a dry region ($\Omega_d$) and a wet region ($\Omega_w$), assuming that water flows only within $\Omega_w$. The governing equation and boundary conditions for a steady-state seepage problem that satisfies Darcy's law are as follows:
\begin{equation}
\nabla(k\nabla h)=0,\quad\text{in}\quad\Omega,
\end{equation}

\begin{equation}
h(x,y,z)=\phi=
\begin{cases}
H_1 & \quad \text{on} \quad S_1 \\
H_2 & \quad \text{on} \quad S_5 \\
z   & \quad \text{on} \quad S_3 \quad \text{and} \quad S_4
\end{cases},
\label{eq:free surface b1}
\end{equation}

\begin{equation}
-k_n\frac{\partial\phi}{\partial n}=\bar{\mathrm{q}}_n=0,\quad\text{on}S_2 \quad \text{and} \quad S_3.
\label{eq:free surface b2}
\end{equation}
where \( S_1 \) and \( S_5 \) are the prescribed hydraulic head boundaries at the upstream and downstream, respectively. \( S_2 \), \( S_3 \), and \( S_4 \) correspond to the prescribed impermeable boundary, the free seepage surface, and the seepage outflow boundary, respectively.

\begin{figure}[H]
  \centering
  \includegraphics[width=0.6\textwidth]{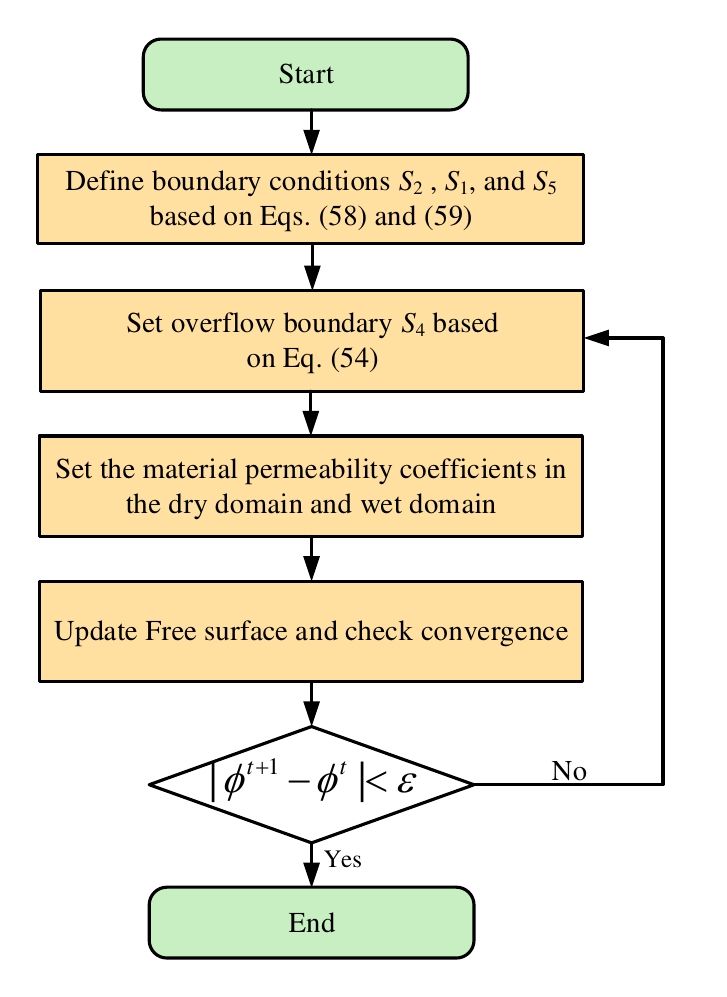}
  \caption{Free surface solution procedure using the PSBFEM.}
  \label{fig:Free_surface_solve}
\end{figure}

In this work, the fixed mesh method~\cite{bathe1979finite} is employed to determine the free surface, as illustrated in Fig.~\ref{fig:Free_surface_solve}. The procedure consists of the following steps:

(1) \textbf{Boundary condition definition}: Specify the hydraulic head boundary conditions, including the impermeable boundary \( S_2 \), and the prescribed head boundaries at the upstream \( S_1 \) and downstream \( S_5 \), as described by Eqs.~(\ref{eq:free surface b1}) and~(\ref{eq:free surface b2}).
    
(2) \textbf{Overflow boundary setup}: Define the overflow boundary \( S_4 \) based on Eq.~(\ref{eq:free surface b1}). The initial configuration ensures that the downstream overflow boundary coincides with the upstream water level.
    
(3) \textbf{Material property assignment}: Assign the material permeability such that the region below the free surface adopts the standard permeability coefficient \( k \), while the region above the free surface is assigned a reduced permeability, typically \( 0.001k \). Export the ABAQUS input file (``.inp'') and perform the seepage analysis using the PSBFEM user-defined element (UEL).
    
(4) \textbf{Free surface update and convergence check}: Extract the seepage field from the ABAQUS output database (``.odb'') using a Python script. Subtract the elevation hydraulic head from the total hydraulic head to obtain the updated coordinates of the free surface and overflow points \( \phi^{t+1} \). If the convergence criterion \( |\phi^{t+1} - \phi^t| < \epsilon \) is satisfied, the computation is considered converged. Otherwise, the procedure returns to Step~(2), the overflow boundary \( S_4 \) is updated accordingly, and the computation is repeated until convergence is achieved.

\section{Numerical examples}
\label{sec:7}
This section presents a series of benchmark problems to evaluate the convergence and accuracy of the proposed framework for three-dimensional seepage analysis. To assess the performance of the PSBFEM, numerical results are compared against those obtained from conventional finite element analysis performed in ABAQUS. All simulations were carried out on a system with an Intel Core i7-4710MQ CPU (2.50~GHz) and 4~GB of RAM. The accuracy of the proposed method is quantified by computing the relative error in hydraulic head, defined as:
\begin{equation}
    \mathbf{e}_{L_2}=\frac{\left\|\mathbf{H}_{num}-\mathbf{H}_{ref}\right\|_{L_2}}{\parallel\mathbf{H}_{ref}\parallel_{L_2}},
\end{equation}
where $\mathbf{H}_{num}$ represents the numerical solution, and $\mathbf{H}_{ref}$ denotes the reference solution.

\subsection{Patch test}
\begin{figure}[H]
  \centering
  \includegraphics[width=0.8\textwidth]{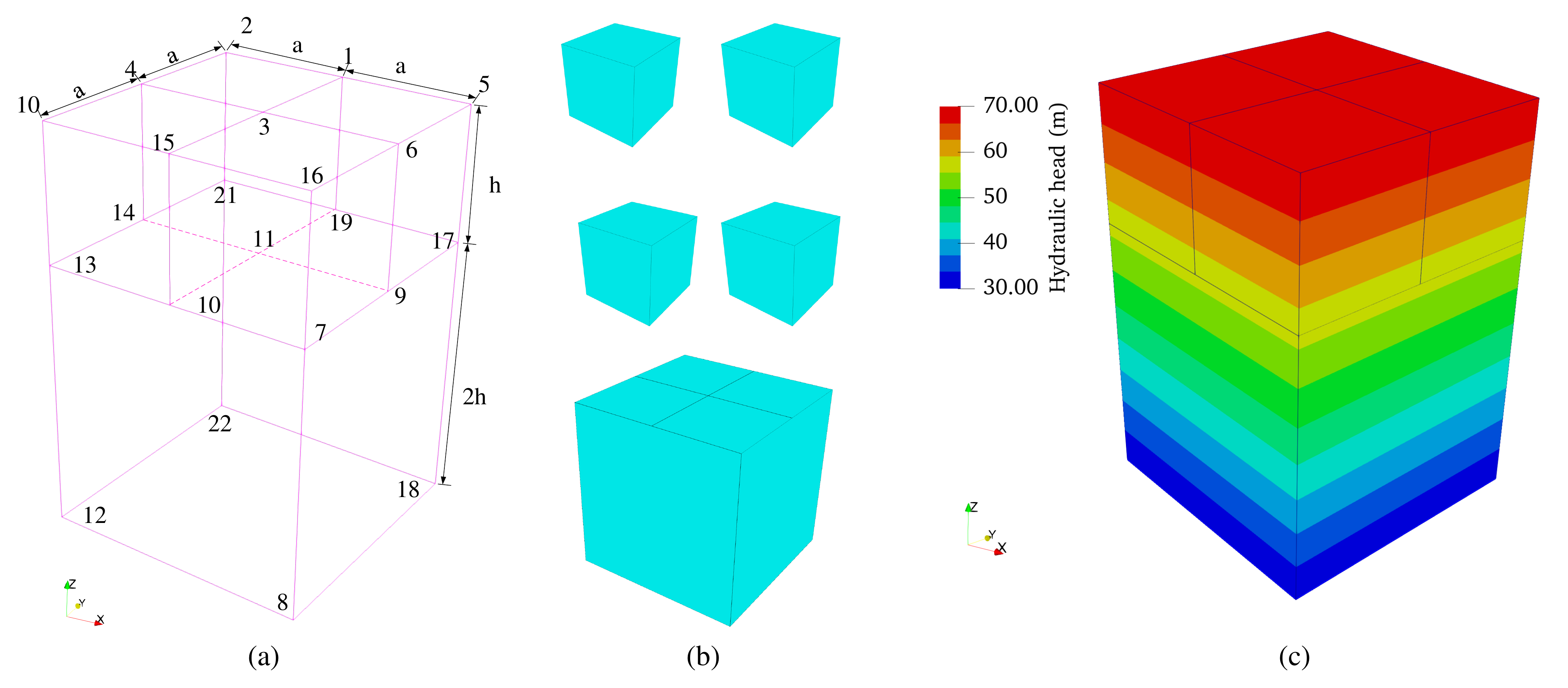}
  \caption{Geometry model and hydraulic head distribution of the patch test; (a) patch element; (b) patch element composition; (c) hydraulic head distribution.}
  \label{fig:patch}
\end{figure}

To examine the convergence behavior of the proposed method, a standard patch test is performed~\cite{ya2021open,zhang2019fast}. The geometry of the test is shown in Fig.~\ref{fig:patch}(a), consisting of a unit-sized quadrangular prism with dimensions \(a = b = h = 0.25\,\text{m}\). As illustrated in Fig.~\ref{fig:patch}(b), the patch comprises four hexahedral elements and a polyhedral element with nine faces. A uniform hydraulic conductivity of \(k = 1 \times 10^{-5}\,\mathrm{m/s}\) is prescribed. Dirichlet boundary conditions are applied with hydraulic heads of 70~m and 30~m at the top (\(z = 3\,\text{m}\)) and bottom (\(z = 0\,\text{m}\)) surfaces, respectively.

The computed hydraulic head distribution is presented in Fig.~\ref{fig:patch}(c). The corresponding relative errors, evaluated at all nodal points against the analytical solution, are listed in Tab.~\ref{tab:patch}. The maximum relative error is found to be \(3.53 \times 10^{-5}\), demonstrating excellent agreement with the analytical solution. These results confirm that the proposed PSBFEM passes the patch test with high accuracy.

\begin{table}[H]
\centering
\caption{Maximum relative error of the nodal hydraulic head.}
\begin{tabular}{@{}lll@{}}
\toprule
PSBFEM (m) & Analytical solution (m)  & Relative error\\ \midrule
 56.6687 & 56.6667 & $3.53 \times 10^{-5}$\\ \bottomrule
\end{tabular}
\label{tab:patch}
\end{table}

\subsection{Steady state seepage analysis}
\subsubsection{Steady state seepage problem in a concrete dam}
\begin{figure}[H]
  \centering
  \includegraphics[width=0.8\textwidth]{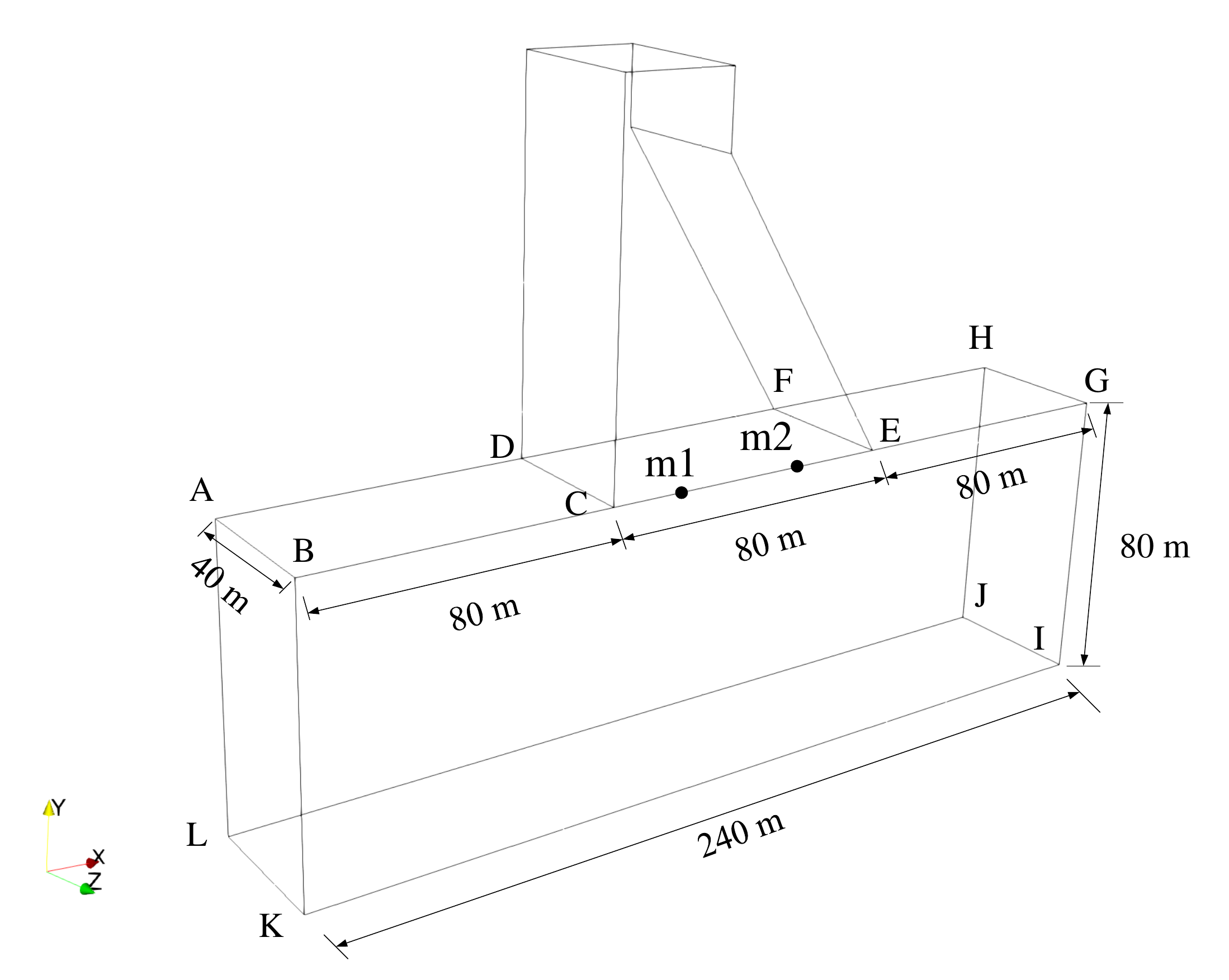}
  \caption{Geometry model of a concrete dam.}
  \label{fig:ex01_geo}
\end{figure}

In this example, a steady-state seepage problem in a concrete dam foundation is analyzed. Since the dam body is assumed to be completely impervious, it is excluded from the numerical model. The geometric configuration is illustrated in Fig.~\ref{fig:ex01_geo}. The upstream boundary \(ABCD\) is prescribed with a hydraulic head of 80~m, and the downstream boundary \(EFGH\) with 20~m. All other boundaries are treated as impermeable. To assess the accuracy of the proposed method, two monitoring points, \(m_1(100, 80, 40)\) and \(m_2(140, 80, 40)\), are selected, as marked in Fig.~\ref{fig:ex01_geo}. The hydraulic conductivity of the dam foundation is set as \(k_x = k_y = k_z = 1 \times 10^{-5}~\mathrm{cm/s}\). For comparison, the ABAQUS uses the C3D8P element, whereas the proposed PSBFEM employs polyhedral elements.

\begin{figure}[H]
  \centering
  \includegraphics[width=0.8\textwidth]{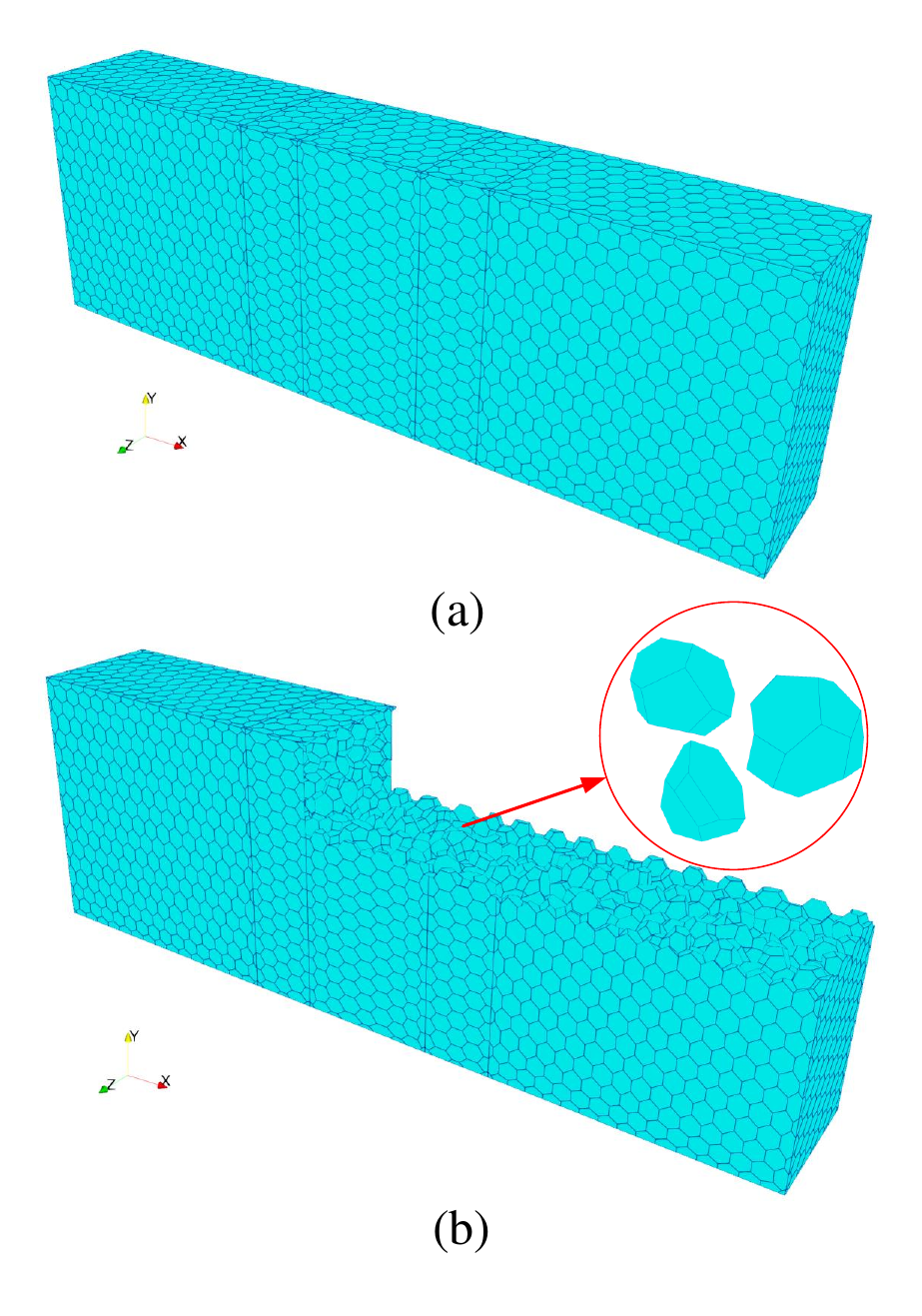}
  \caption{Meshes of concrete dam foundation. (a) 5 m mesh; (b) internal mesh of 5 m mesh size.}
  \label{fig:ex01_P_mesh}
\end{figure}

\begin{figure}[H]
  \centering
  \includegraphics[width=1.0\textwidth]{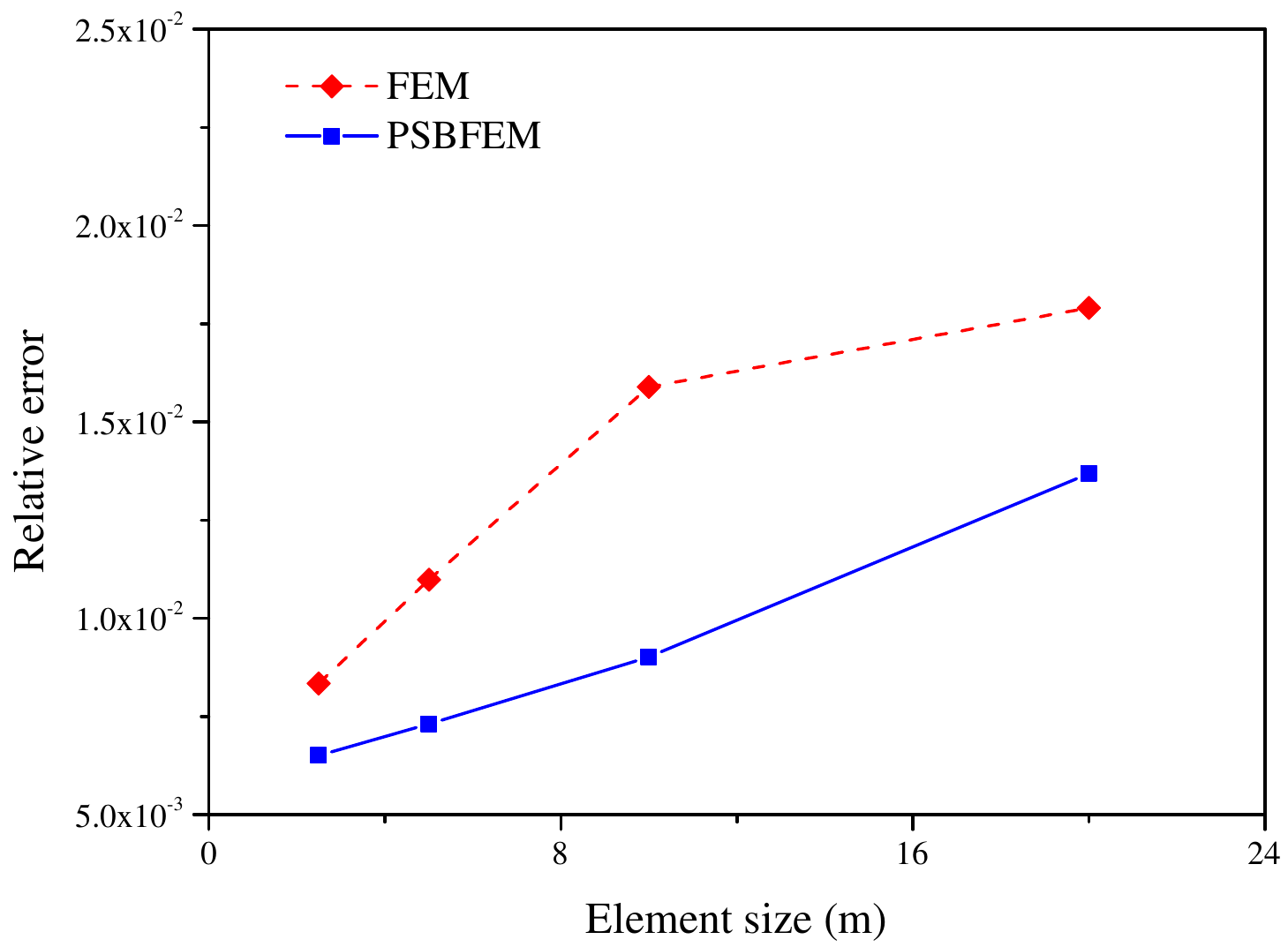}
  \caption{Comparison of the convergence rates in the hydraulic head for monitoring point.}
  \label{fig:ex01_error}
\end{figure}

\begin{figure}[H]
  \centering
  \includegraphics[width=1.0\textwidth]{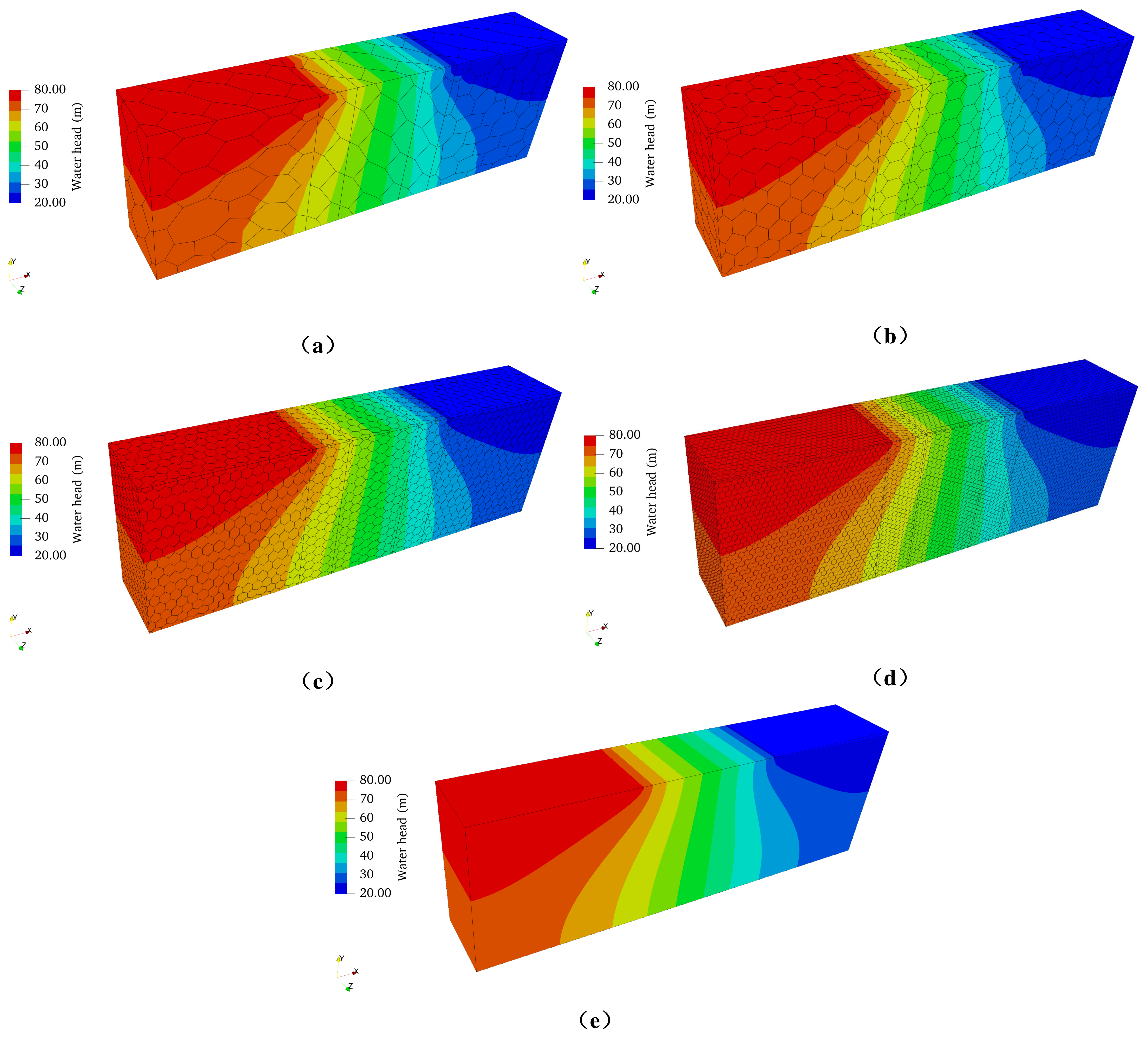}
  \caption{Distribution of hydraulic head for the polyhedral mesh; (a) 20.0 m mesh size; (b) 10.0 m mesh size; (c) 5.0 m mesh size; (d) 2.5 m mesh size; (e) reference solution.}
  \label{fig:ex01_poly_contour}
\end{figure}

A convergence study is conducted through mesh \(h\)-refinement, with mesh sizes of 20~m, 10~m, 5~m, and 2.5~m. The 5~m mesh model is shown in Fig.~\ref{fig:ex01_P_mesh}(a), and Fig.~\ref{fig:ex01_P_mesh}(b) presents a cross-sectional view with representative polyhedral elements. The relative errors of the hydraulic head at the monitoring points are plotted in Fig.~\ref{fig:ex01_error}, revealing that both the PSBFEM and the FEM demonstrate satisfactory convergence. Table~\ref{tab:ex01_t1} lists the computed hydraulic head values at the monitoring points. The relative errors for the FEM and PSBFEM are \(1.5 \times 10^{-2}\) and \(0.9 \times 10^{-2}\), respectively, indicating that the PSBFEM achieves higher accuracy under identical mesh resolutions. Furthermore, Fig.~\ref{fig:ex01_poly_contour} shows that the hydraulic head contours obtained by the PSBFEM closely match the reference solution for the refined mesh.

\begin{figure}[H]
  \centering
  \includegraphics[width=1.0\textwidth]{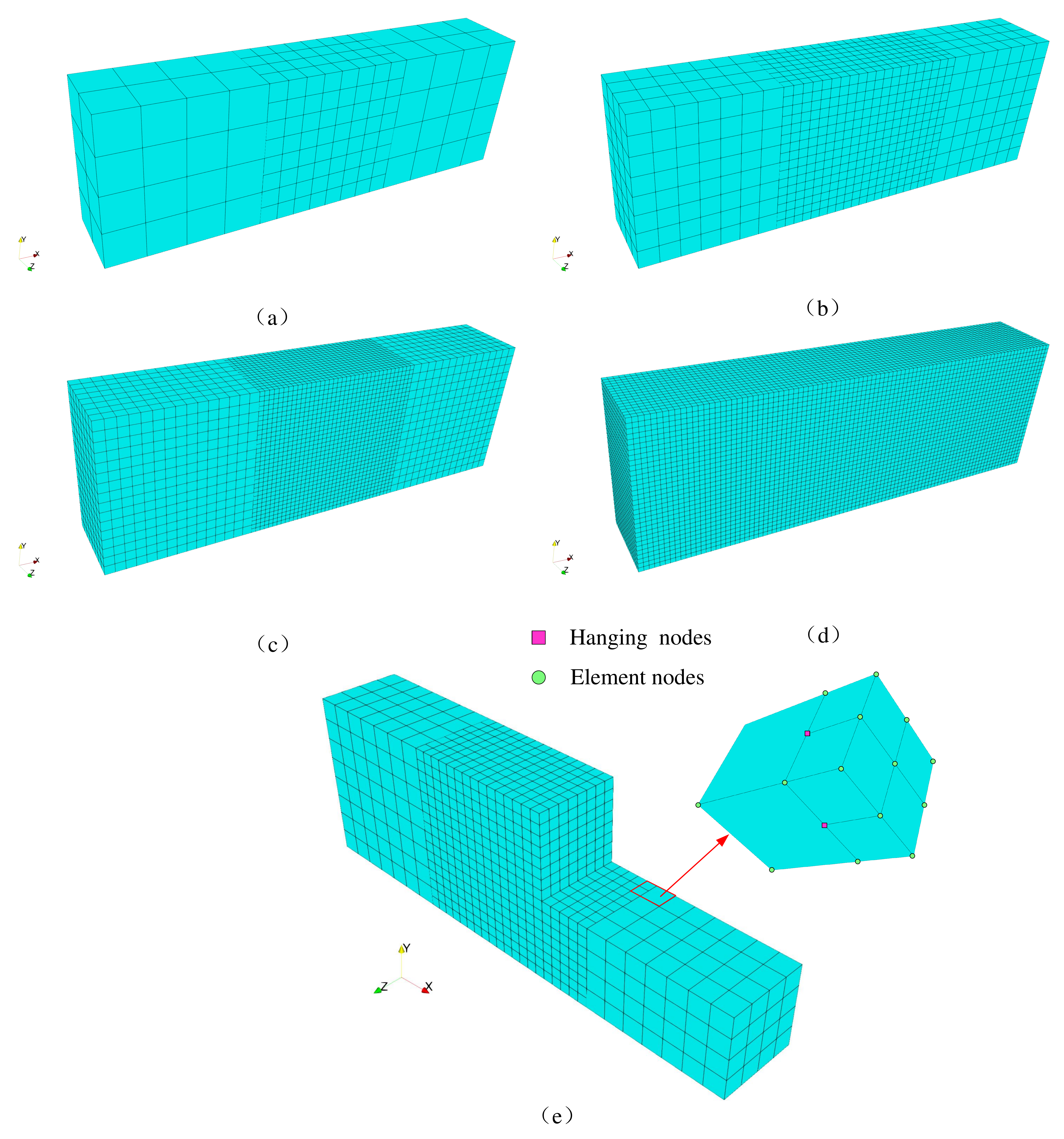}
  \caption{Distribution of hydraulic head for the polyhedral mesh; (a) Refining mesh level 1; (b) Refining mesh level 2; (c) Refining mesh level 3; (d) fine mesh; (e) internal mesh of refining mesh level 3.}
  \label{fig:ex01_octree}
\end{figure}

\begin{figure}[H]
  \centering
  \includegraphics[width=1.0\textwidth]{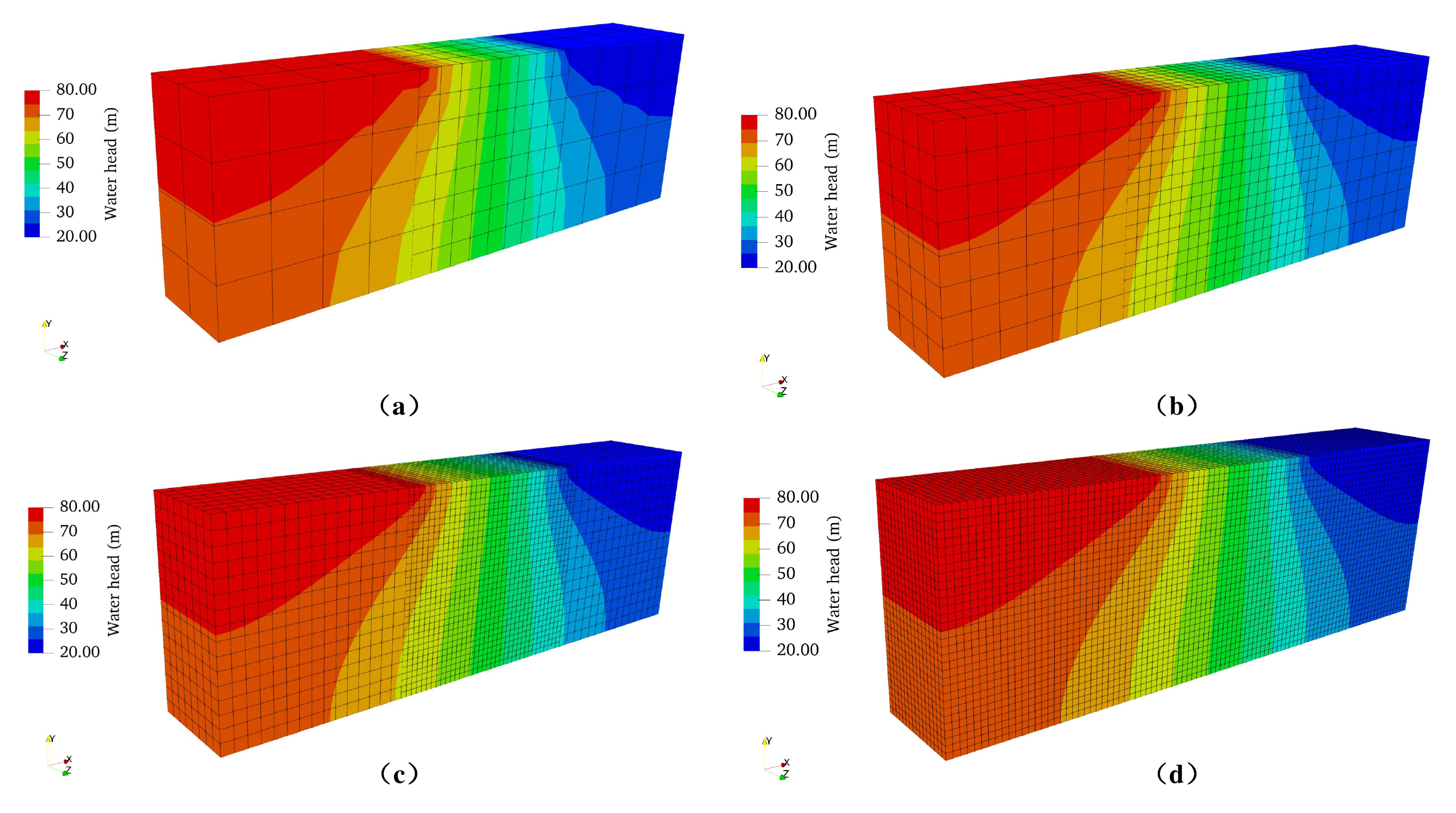}
  \caption{Distribution of hydraulic head for the octree mesh; (a) refining mesh 1; (b) refining mesh 2; (c) refining mesh 3; (d) fine mesh.}
  \label{fig:ex01_octree_contour}
\end{figure}

To better capture the seepage behavior beneath the lower part of the dam, local mesh refinement is applied to the dam foundation using an octree-based strategy, as shown in Figs.~\ref{fig:ex01_octree}(a)–(c). A cross-sectional view of the level-3 refined mesh is illustrated in Fig.~\ref{fig:ex01_octree}(e). The mesh characteristics and corresponding relative errors are summarized in Table~\ref{tab:ex01_t2}. The relative error for level-3 refinement is comparable to that of the globally refined mesh, yet with significantly reduced computational cost. This demonstrates that the octree-based local refinement effectively balances accuracy and efficiency. Moreover, the hydraulic head distribution obtained using the PSBFEM with the octree mesh shows excellent agreement with the reference solution, as depicted in Fig.~\ref{fig:ex01_octree_contour}.

\begin{table}[H]
    \centering
    \caption{Comparison of the hydraulic head using the PSBFEM and FEM (element size: 10 m).}
    \label{tab:comparison}
    \begin{tabular}{cccc}
        \hline
        Methods & Monitor point m1  & Monitor point m2  & Relative error \\ 
        \hline
        Analytical solutions \cite{yangNovel2022} & 60.00 & 40.00 & -- \\
        FEM & 60.7642 & 39.2387 & 1.5$\times10^{-2}$ \\	
        PSBFEM & 60.5119 & 39.6203 & 0.9$\times10^{-2}$ \\	
        \hline
        \label{tab:ex01_t1}
    \end{tabular}
\end{table}

\begin{table}[H]
\centering
\caption{Mesh characteristics and relative errors for the locally refined meshes.}
\begin{tabular}{cccccc}
\toprule
Mesh type  & Elements & Nodes & Faces &  Relative error & CPU time (s)\\ 
\midrule
Refining mesh level 1            & 320       & 525       & 1968       & 1.59$\times10^{-2}$  & 0.70  \\ 
Refining mesh level 2     & 3008      & 3843      & 18240       &  1.11$\times10^{-2}$ &    3.40\\ 
Refining mesh level 3             & 20256      & 29001     & 128688     & 8.47$\times10^{-3}$   &  36.90  \\ 
Fine mesh      & 21264      & 24671      & 129072     & 8.43$\times10^{-3}$   & 41.70   \\ 
\bottomrule 
\label{tab:ex01_t2}
\end{tabular}
\end{table}

\subsubsection{Steady state seepage analysis for permeable materials}
\begin{figure}[H]
  \centering
  \includegraphics[width=1.0\textwidth]{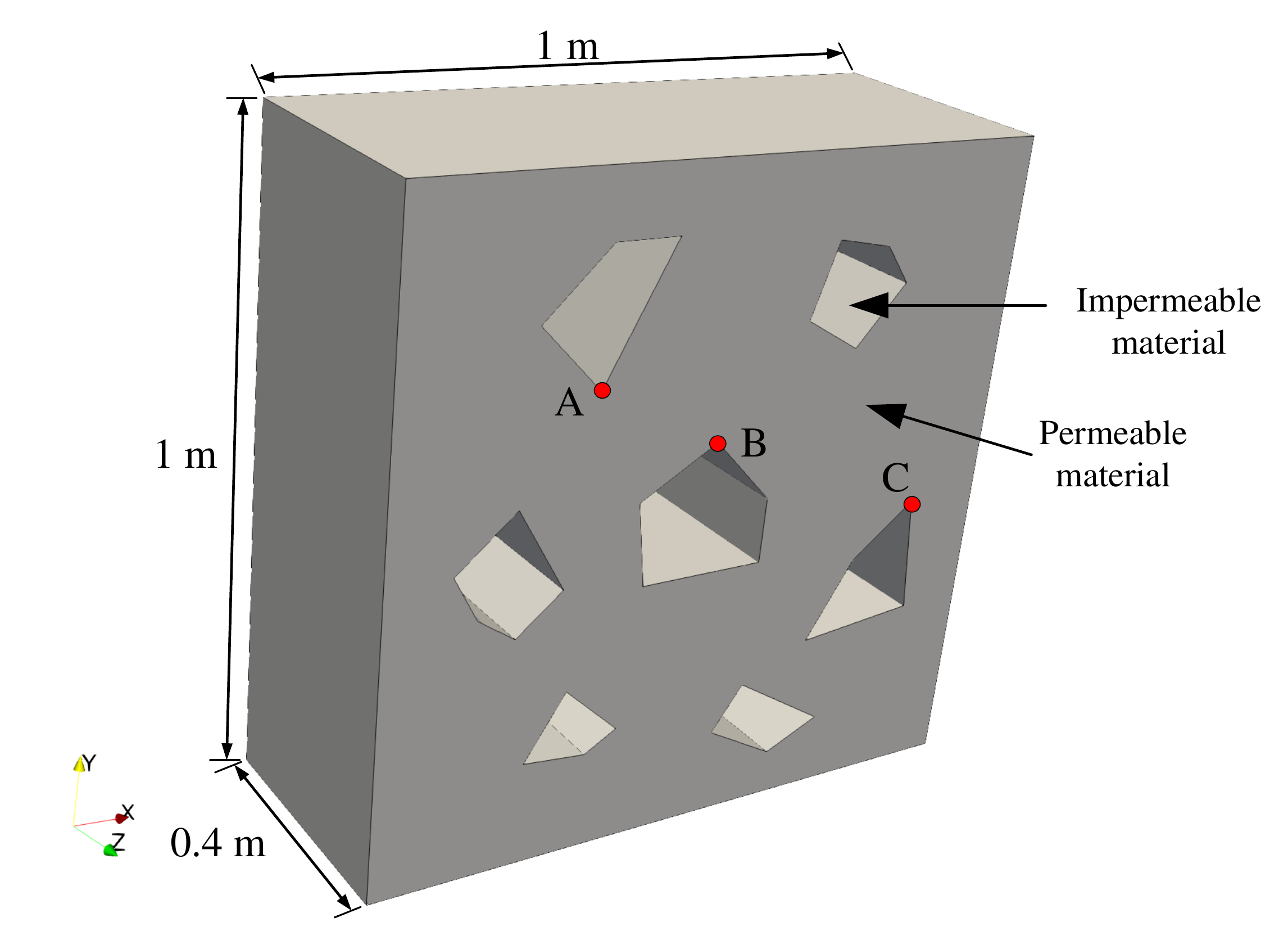}
  \caption{Permeable material’s geometric model.}
  \label{fig:ex02_geo}
\end{figure}

\begin{figure}[H]
  \centering
  \includegraphics[width=1.0\textwidth]{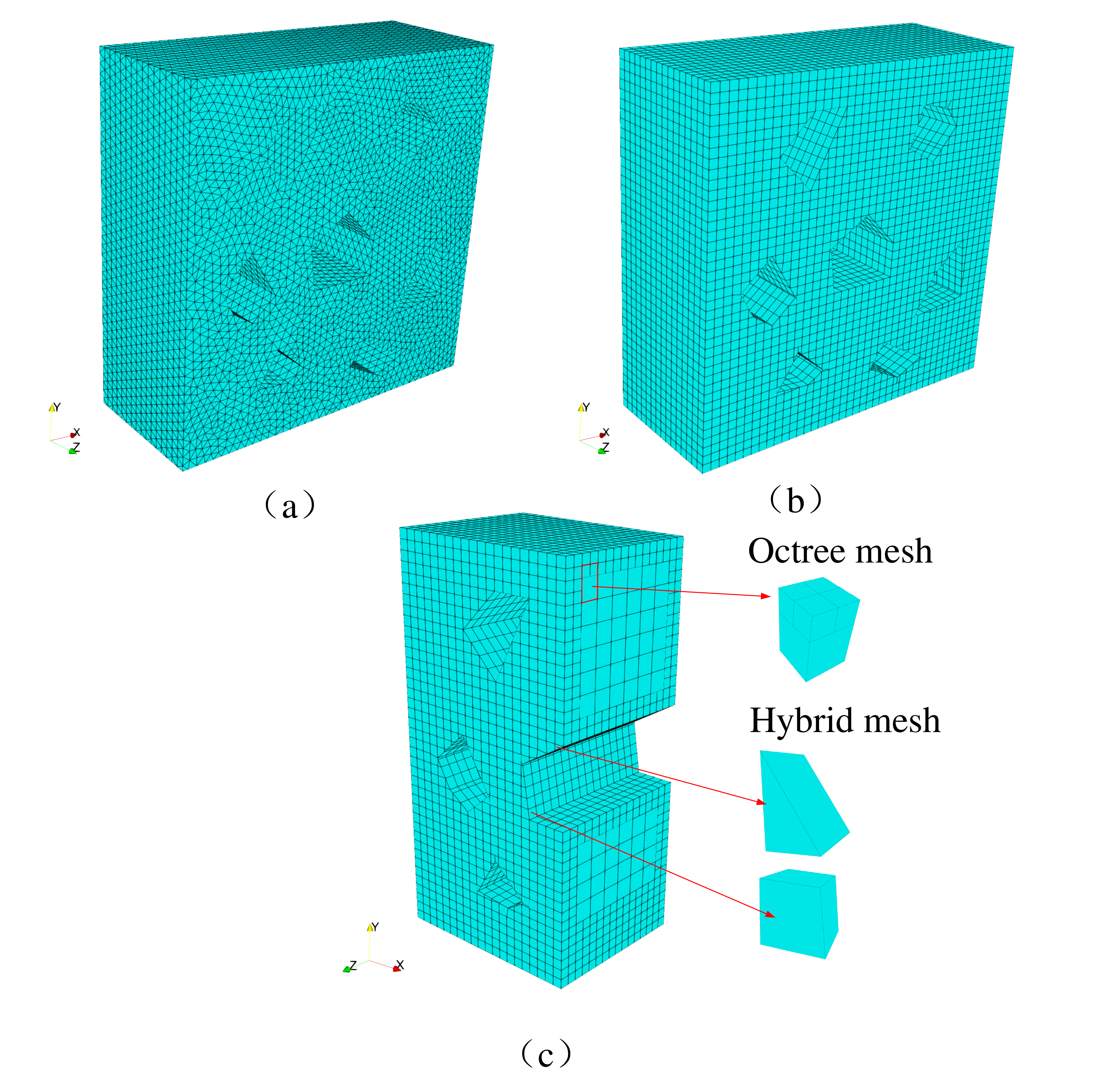}
  \caption{Permeable material’s mesh model; (a) tetrahedral mesh; (b) hybrid octree mesh. (c) internal mesh of
hybrid octree mesh.}
  \label{fig:ex02_mesh}
\end{figure}

To demonstrate the adaptability of the proposed PSBFEM integrated with a hybrid octree mesh, a steady-state seepage analysis is performed on a permeable medium containing an embedded impermeable inclusion, as illustrated in Fig.~\ref{fig:ex02_geo}. The domain has a length and width of 1~m, and the permeability is isotropic with $k_x = k_y = k_z = 1 \times 10^{-5}$~cm/s. The hybrid octree mesh, shown in Fig.~\ref{fig:ex02_mesh}(b), effectively captures the geometric features of the domain. In regions of geometric regularity, it generates structured hexahedral elements, as depicted in Fig.~\ref{fig:ex02_mesh}(c), while irregular polyhedral elements are employed near boundaries and transition zones. For comparison, the conventional FEM is applied using an unstructured tetrahedral mesh, as shown in Fig.~\ref{fig:ex02_mesh}(a). For both mesh types, the element size is set to 0.025~m.

\begin{figure}[H]
  \centering
  \includegraphics[width=1.0\textwidth]{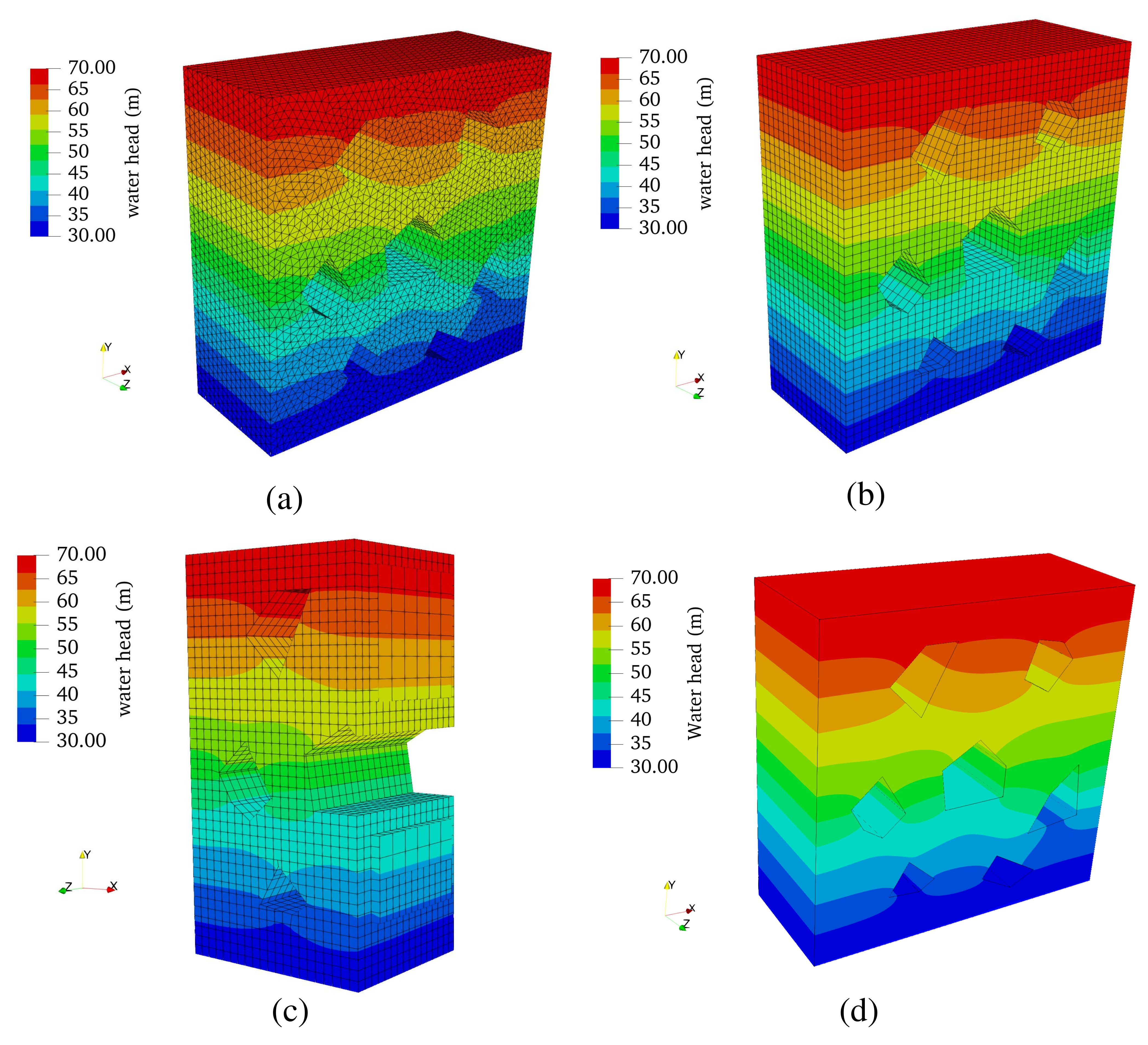}
  \caption{Distribution of hydraulic head for the permeable material; (a) the solution of FEM; (b) the solution of PSBFEM; (c) localized contour of selected octree elements from (b); (d) the reference solution.}
  \label{fig:ex02_contour}
\end{figure}

The hydraulic head distributions obtained from PSBFEM and FEM are presented in Fig.~\ref{fig:ex02_contour}, showing close agreement with the reference solution. As summarized in Table~\ref{tab:ex02_t1}, the relative errors at selected monitoring points obtained by PSBFEM are consistently on the order of $10^{-3}$, whereas those obtained by FEM are on the order of $10^{-2}$. These results indicate that PSBFEM achieves superior numerical accuracy compared to FEM in this problem.

\begin{table}[H]
\centering
\caption{Relative errors for the monitoring points.}
\label{tab:ex02_t1}
\resizebox{\textwidth}{!}{%
\begin{tabular}{cccccc}
\toprule
Monitoring point & PSBFEM (m) & FEM (m) & Reference solution & $e_{L2}$ of PSBFEM & $e_{L2}$ of FEM \\ 
\midrule
A & 57.3925 & 57.7907 & 57.1854    & $3.62\times10^{-3}$ & $1.06\times10^{-2}$ \\
B & 56.9803 & 57.3528 & 56.7711      & $3.68\times10^{-3}$ & $1.03\times10^{-2}$ \\
C & 50.0317 & 50.2813 & 49.6992  & $6.69\times10^{-3}$ & $1.17\times10^{-2}$ \\
\bottomrule
\end{tabular}%
}
\end{table}

\subsection{Transient seepage analysis}
\subsubsection{Transient seepage analysis of dam foundation with irregular geometry}
\begin{figure}[H]
  \centering
  \includegraphics[width=0.8\textwidth]{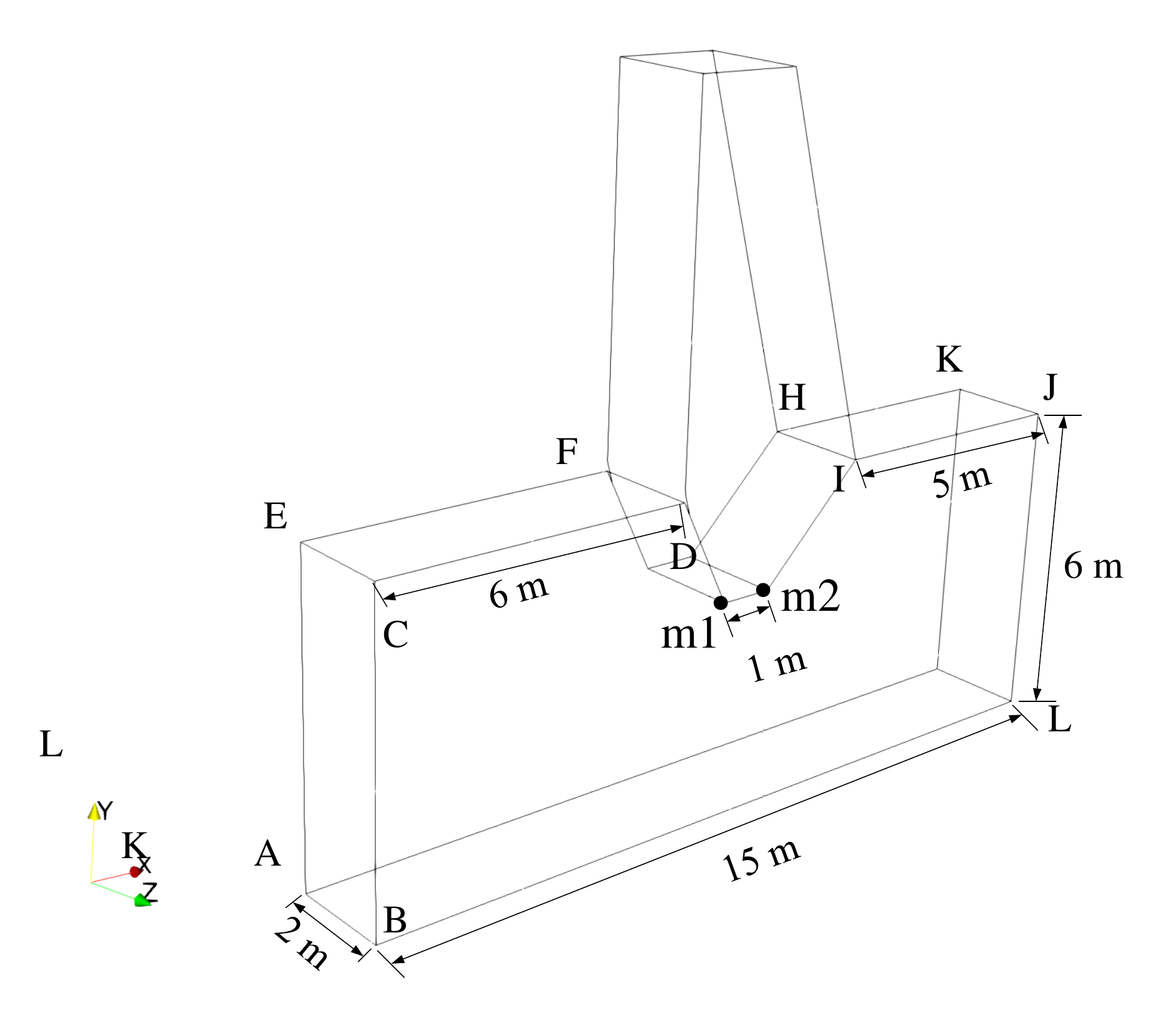}
  \caption{Geometry of a concrete dam with irregular geometry.}
  \label{fig:ex03_geo}
\end{figure}

\begin{figure}[H]
  \centering
  \includegraphics[width=1.0\textwidth]{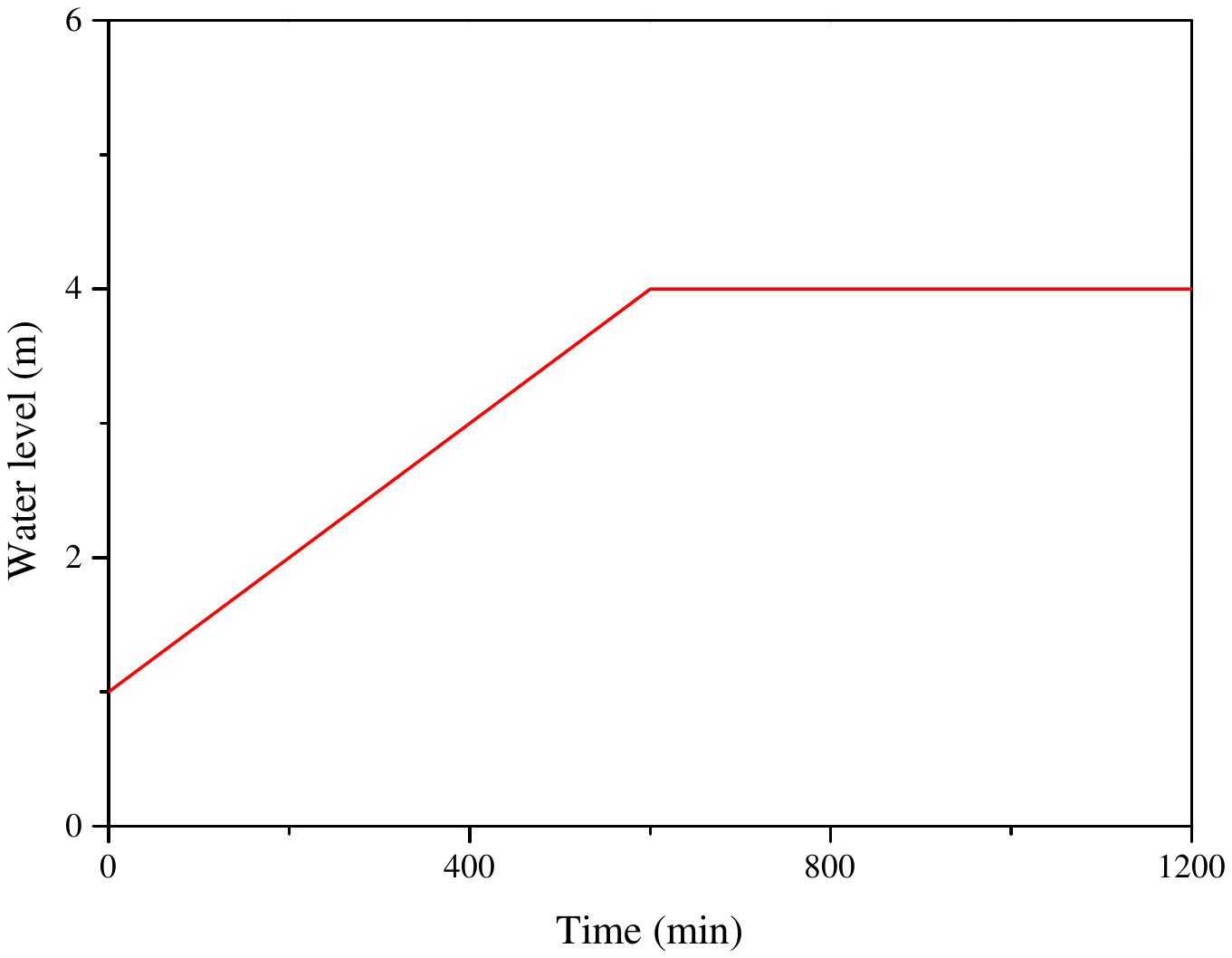}
  \caption{Hydraulic head boundary condition in the upstream.}
  \label{fig:ex03_waterLevel}
\end{figure}

This example examines a concrete dam with an irregular geometry, as depicted in Fig.~\ref{fig:ex03_geo}. To evaluate the transient hydraulic response within the dam foundation, two monitoring points, $m_1$ and $m_2$, are selected. The storage coefficient is set to $S_s = 1.0 \times 10^{-3}\,\mathrm{m}^{-1}$, and the permeability coefficient is $k = 1.02 \times 10^{-3}\,\mathrm{m/min}$. Initially, the upstream and downstream water levels are 1.0~m and 0.0~m, respectively. The upstream level increases gradually to 4.0~m over time, as shown in Fig.~\ref{fig:ex03_waterLevel}. The transient simulation is conducted over a total duration of 1200 minutes using a constant time step of 10 minutes.

\begin{figure}[H]
  \centering
  \includegraphics[width=1.0\textwidth]{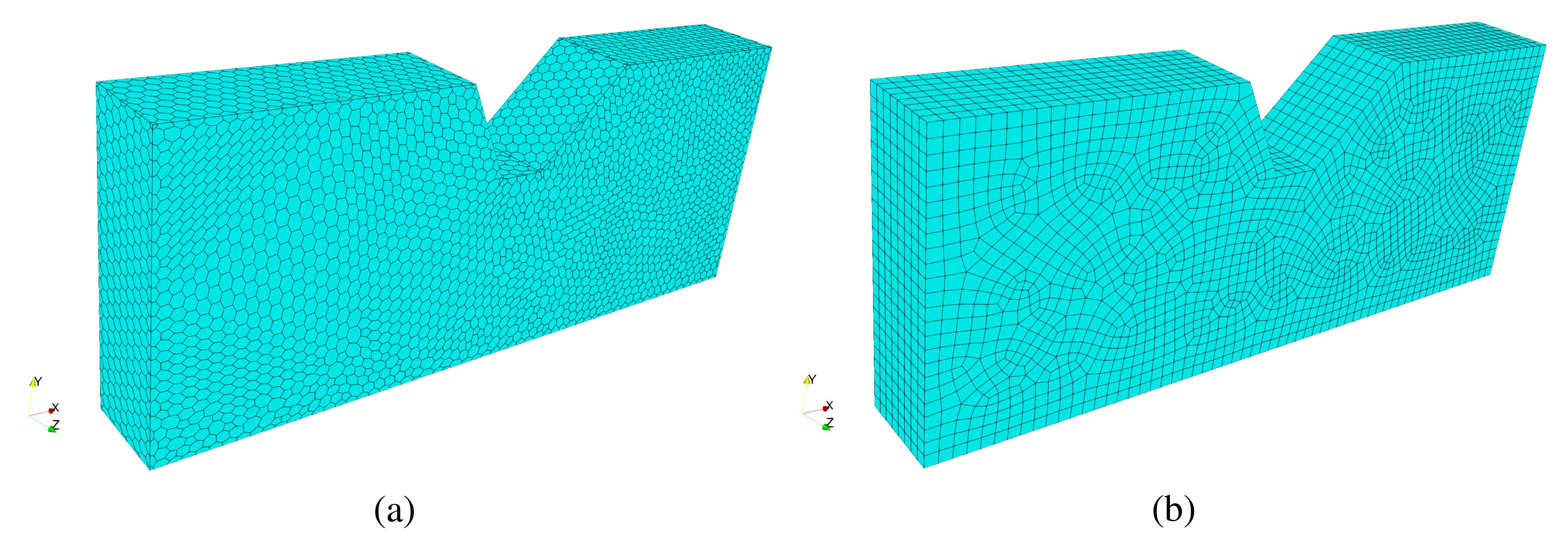}
  \caption{The Mesh of dam foundation with irregular geometry; (a) polyhedral mesh; (b) hexahedronral mesh.}
  \label{fig:ex03_mesh}
\end{figure}

\begin{figure}[H]
  \centering
  \includegraphics[width=0.8\textwidth]{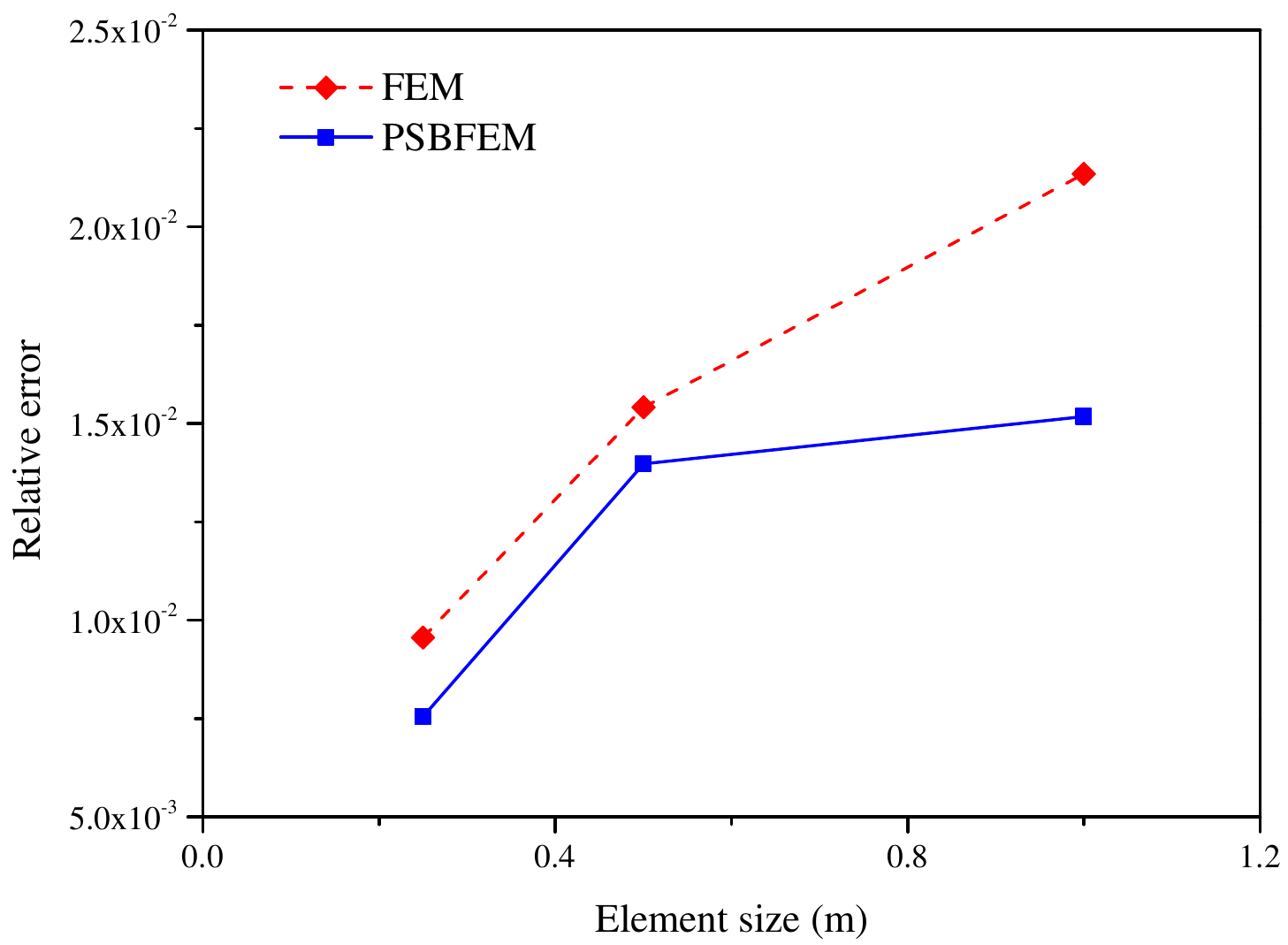}
  \caption{Comparison of the hydraulic head history at the monitoring point.}
  \label{fig:ex03_error}
\end{figure}

\begin{figure}[H]
  \centering
  \includegraphics[width=0.8\textwidth]{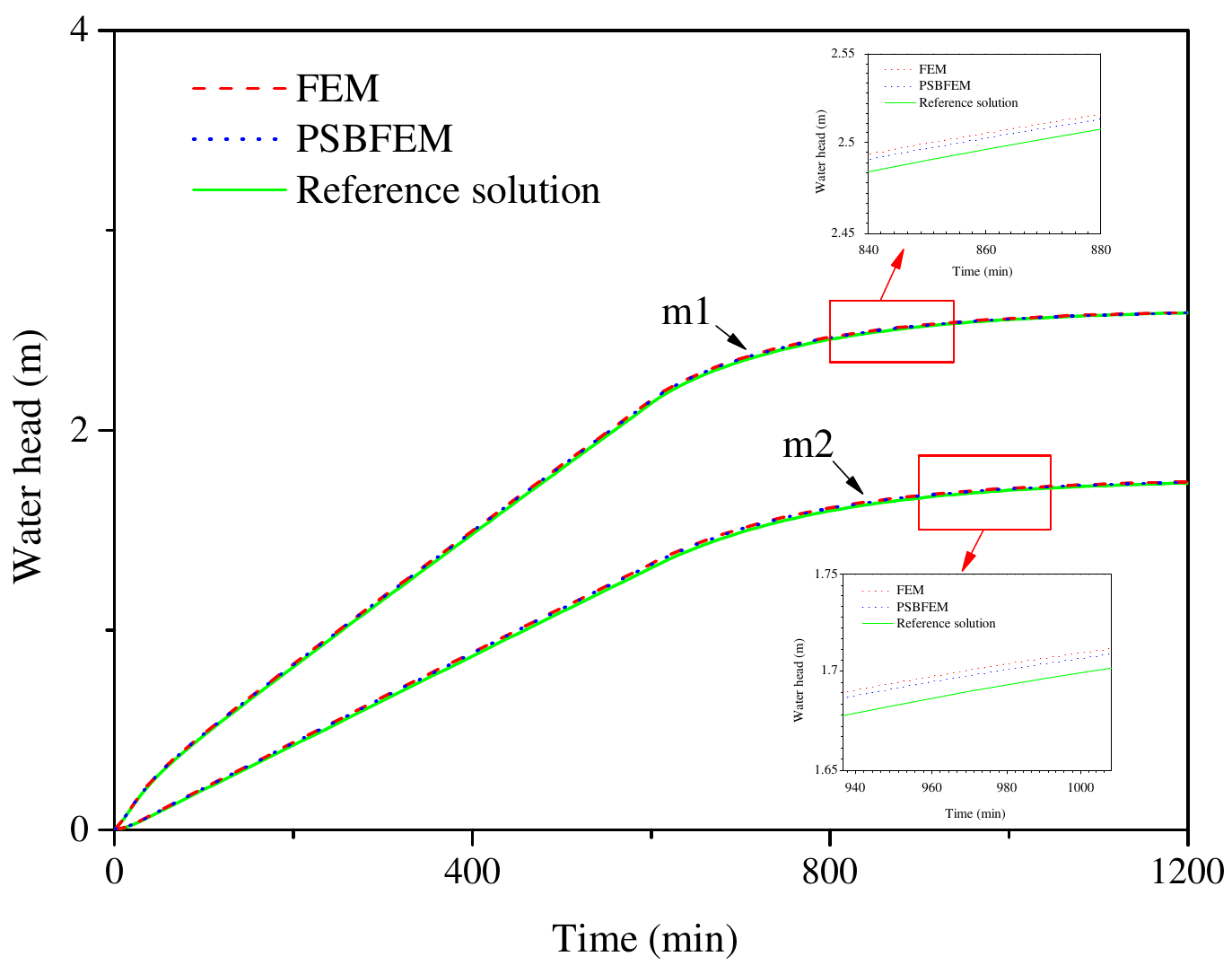}
  \caption{Comparison of the hydraulic head history at the monitoring point (element size 0.25 m).}
  \label{fig:ex03_his}
\end{figure}

To assess spatial convergence, three mesh resolutions with element sizes of 1.0~m, 0.5~m, and 0.25~m are employed. The finest mesh (0.25~m) is illustrated in Fig.~\ref{fig:ex03_mesh}. Fig.~\ref{fig:ex03_his} presents the temporal variation of the hydraulic head at the monitoring points, where the PSBFEM results show excellent agreement with the reference solution. The convergence behavior with mesh refinement is depicted in Fig.~\ref{fig:ex03_error}, indicating that PSBFEM exhibits superior accuracy and a faster convergence rate compared to the conventional FEM.

\begin{figure}[H]
  \centering
  \includegraphics[width=1.0\textwidth]{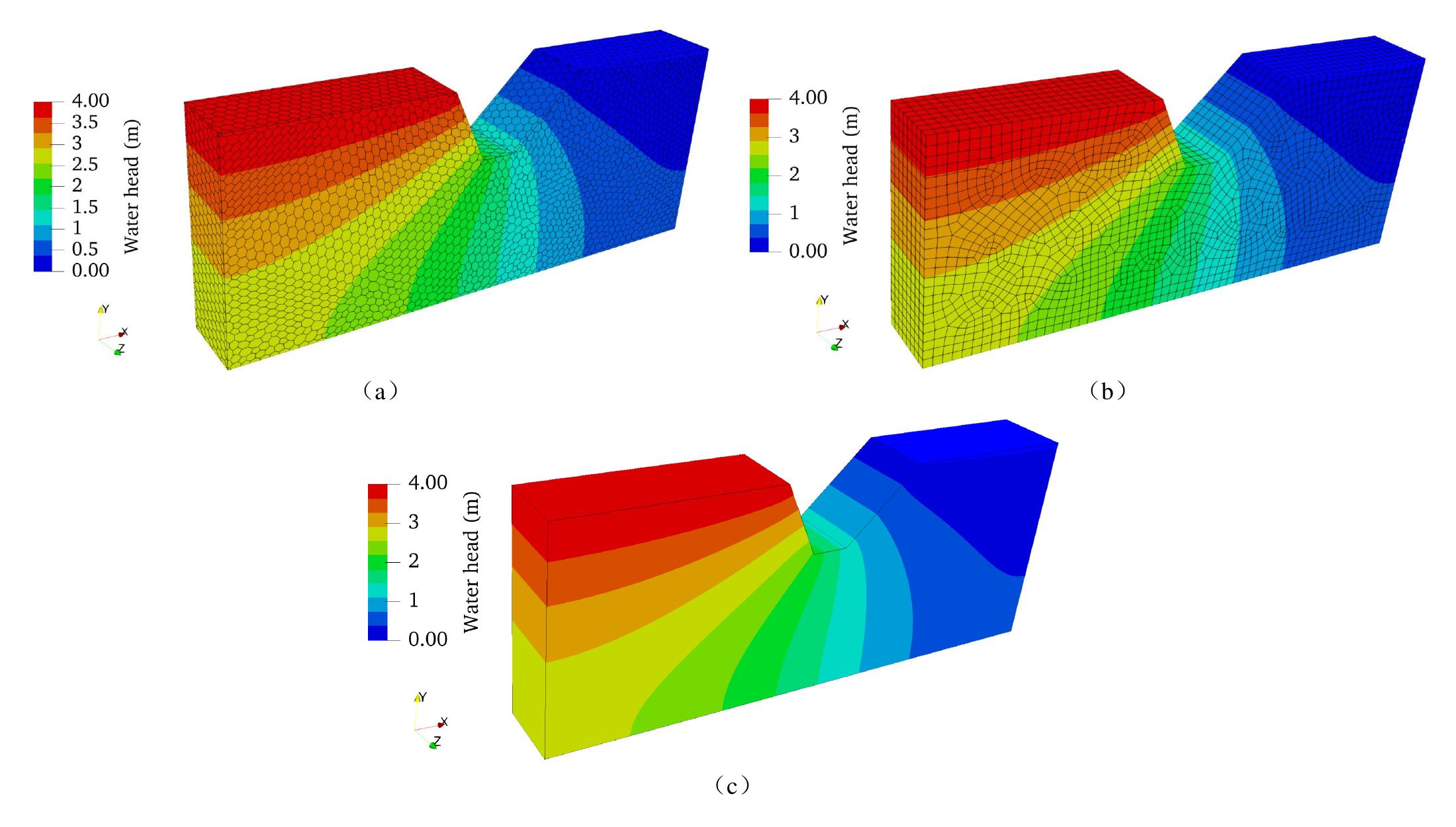}
  \caption{hydraulic head distribution of am foundation with irregular geometry for the transient seepage problem: (a) PSBFEM; (b) FEM; (c) reference solution (element size 0.25 m).}
  \label{fig:ex03_contour}
\end{figure}

Furthermore, Fig.~\ref{fig:ex03_contour} compares the hydraulic head distributions obtained from FEM and PSBFEM. Both methods produce results that are consistent with the reference data, with PSBFEM demonstrating improved accuracy in resolving the hydraulic gradient under refined mesh conditions.

\subsubsection{Transient seepage analysis of complex geometry}
\begin{figure}[H]
  \centering
  \includegraphics[width=1.0\textwidth]{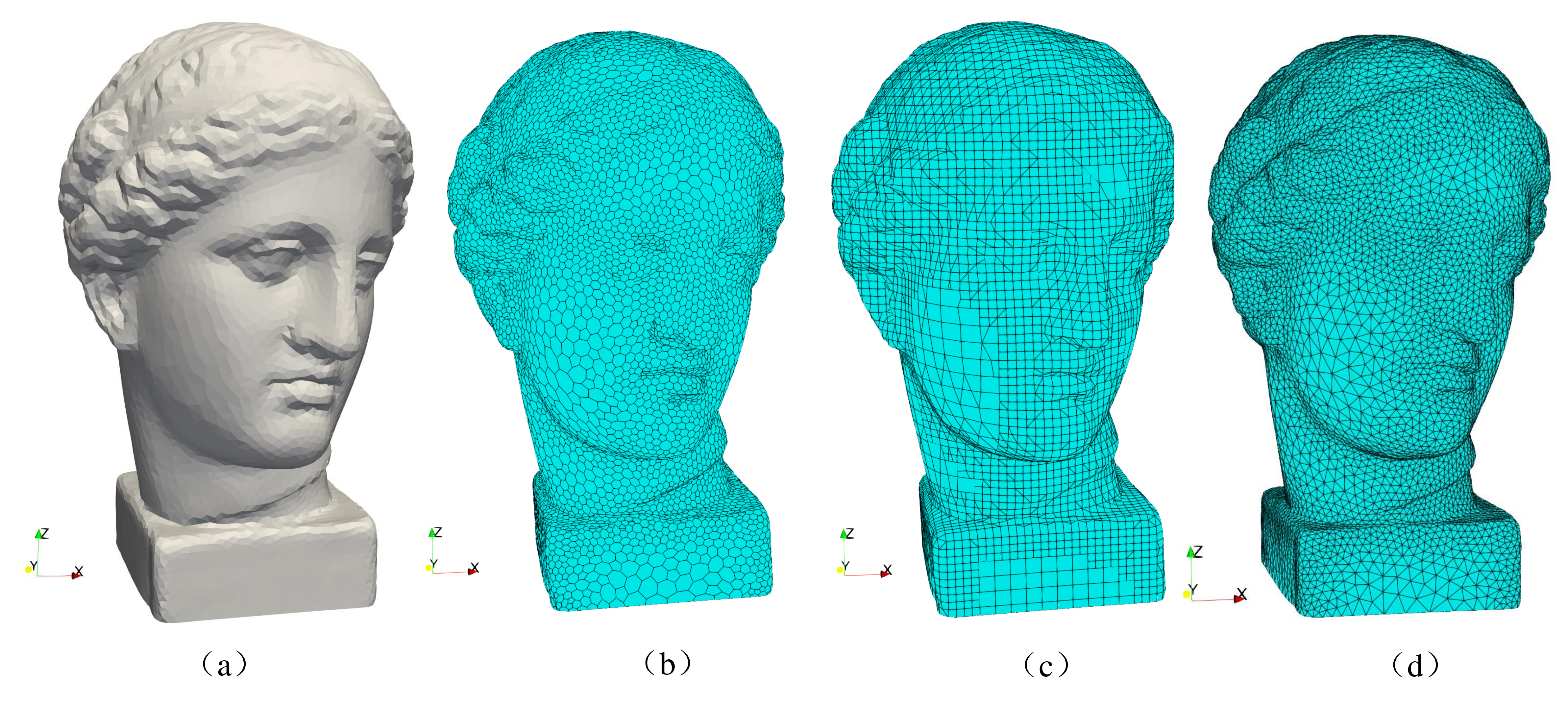}
   \caption{Geometry and mesh model of Venus; (a) the geometry model; (b) the polyhedral mesh; (c) hybrid octree mesh; (d) the tetrahedral mesh.}
  \label{fig:ex04_geo_mesh}
\end{figure}

\begin{figure}[H]
  \centering
  \includegraphics[width=1.0\textwidth]{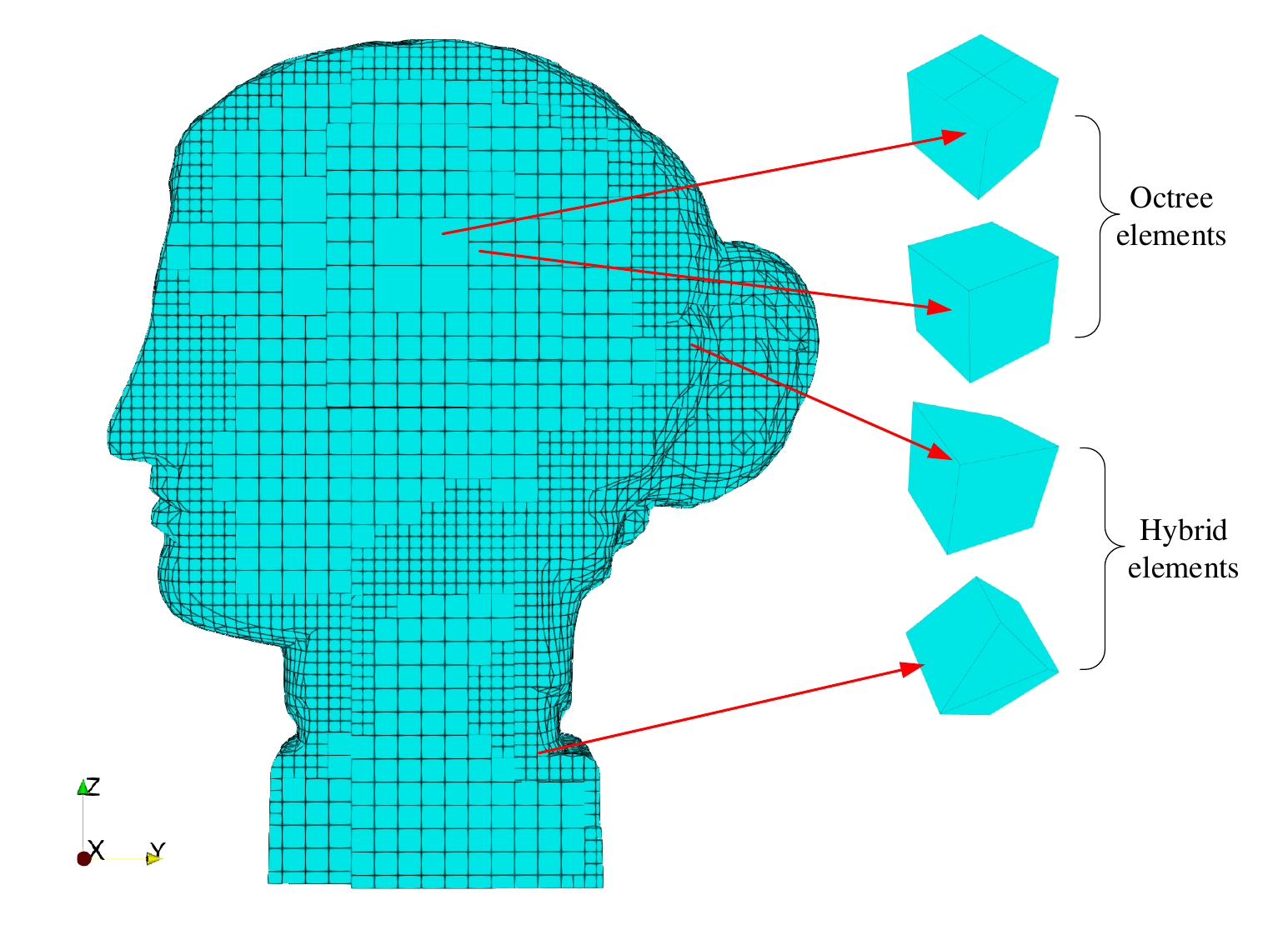}
   \caption{Cross-section view of the hybrid octree mesh for Venus.}
  \label{fig:ex04_hybridOctreeSection}
\end{figure}

In this example, a transient seepage analysis is performed on a model with complex geometry. The geometry is based on the Venus model~\cite{venus_of_arles_3axis}, provided in STL format, and is discretized using polyhedral and hybrid octree meshes, as illustrated in Fig.~\ref{fig:ex04_geo_mesh}(b) and (c). For comparison, the model is also meshed using an unstructured tetrahedral discretization, as shown in Fig.~\ref{fig:ex04_geo_mesh}(d). A cross-sectional view of the hybrid octree mesh is presented in Fig.~\ref{fig:ex04_hybridOctreeSection}, where the interior is filled with regular octree elements and the outer surface is composed of hybrid polyhedral elements generated through a surface-cutting operation. The composition of the hybrid octree mesh is summarized in Tab.~\ref{tab:Composition of hybrid octree mesh}, indicating that 75.9\% of the elements are regular octree cubes and the remaining 24.1\% are hybrid elements. This high proportion of regular elements contributes to the overall mesh quality and numerical accuracy.

\begin{figure}[H]
  \centering
  \includegraphics[width=1.0\textwidth]{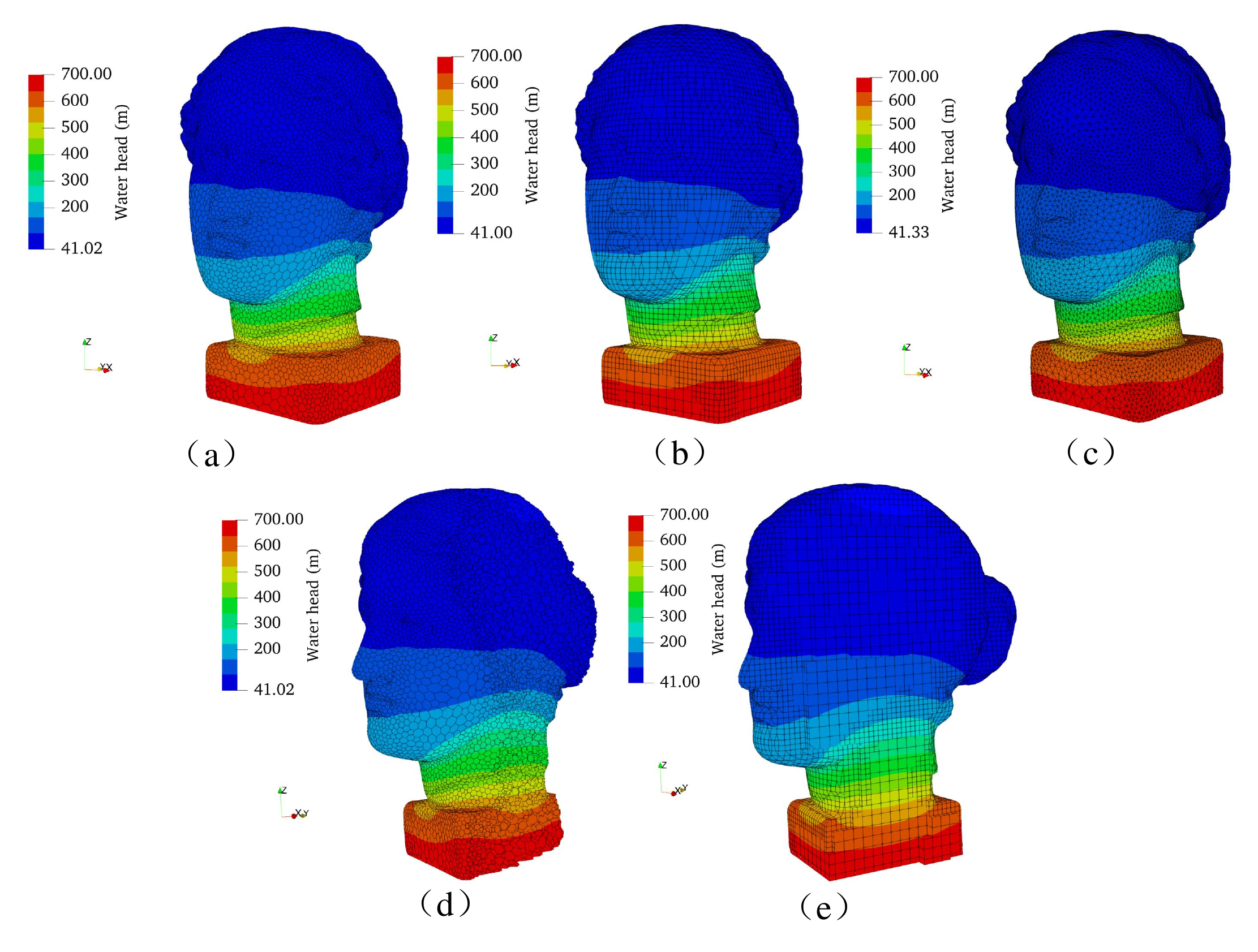}
   \caption{Hydraulic head distribution of Venus at the 300 s. (a) hydraulic head contour with polyhedral elements; (b) hydraulic head contour with hybrid octree elements; (c) reference solution; (d) localized contour of selected polyhedral element; (e) localized contour of selected hybrid octree element.}
  \label{fig:ex04_contour}
\end{figure}

The permeability coefficient is set to $k = 1.7 \times 10^{-4}$~m/s. A hydraulic head of 3~m is applied to the model at the initial time, while a fixed head of 70~m is imposed at the bottom boundary. The resulting hydraulic head distributions obtained from the three different mesh types are shown in Fig.~\ref{fig:ex04_contour}, demonstrating good agreement across all cases and verifying the robustness of the proposed method in handling complex geometries.

\begin{table}[htbp]
\centering
\caption{Composition of hybrid octree mesh for the Venus.}
\begin{tabular}{lccc}
\toprule
Element type & octree element & hybrid element\\
\midrule
Number of elements &22314	&7095 \\
Proportion &75.9\%	&24.1\% \\
\bottomrule
\end{tabular}
\label{tab:Composition of hybrid octree mesh}
\end{table}

\begin{table}[htbp]
\centering
\caption{Mesh characteristics of three elements type.}
\begin{tabular}{lcccc}
\toprule
Element type & Nodes & Elements & Surfaces & CPU time (s) \\
\midrule
Polyhedral mesh & 174104 & 30214 & 395465 &  1798.5\\
Hybrid octree mesh &38077& 29409&184392 & 733.30\\
Tetrahedral mesh & 29358 & 145345 & 581380 &  565.20\\
\bottomrule
\end{tabular}
\label{tab:Comparison of meshes}
\end{table}

\subsection{Free surface seepage problems}
\subsubsection{Homogeneous rectangular dam}
\begin{figure}[H]
  \centering
  \includegraphics[width=1.0\textwidth]{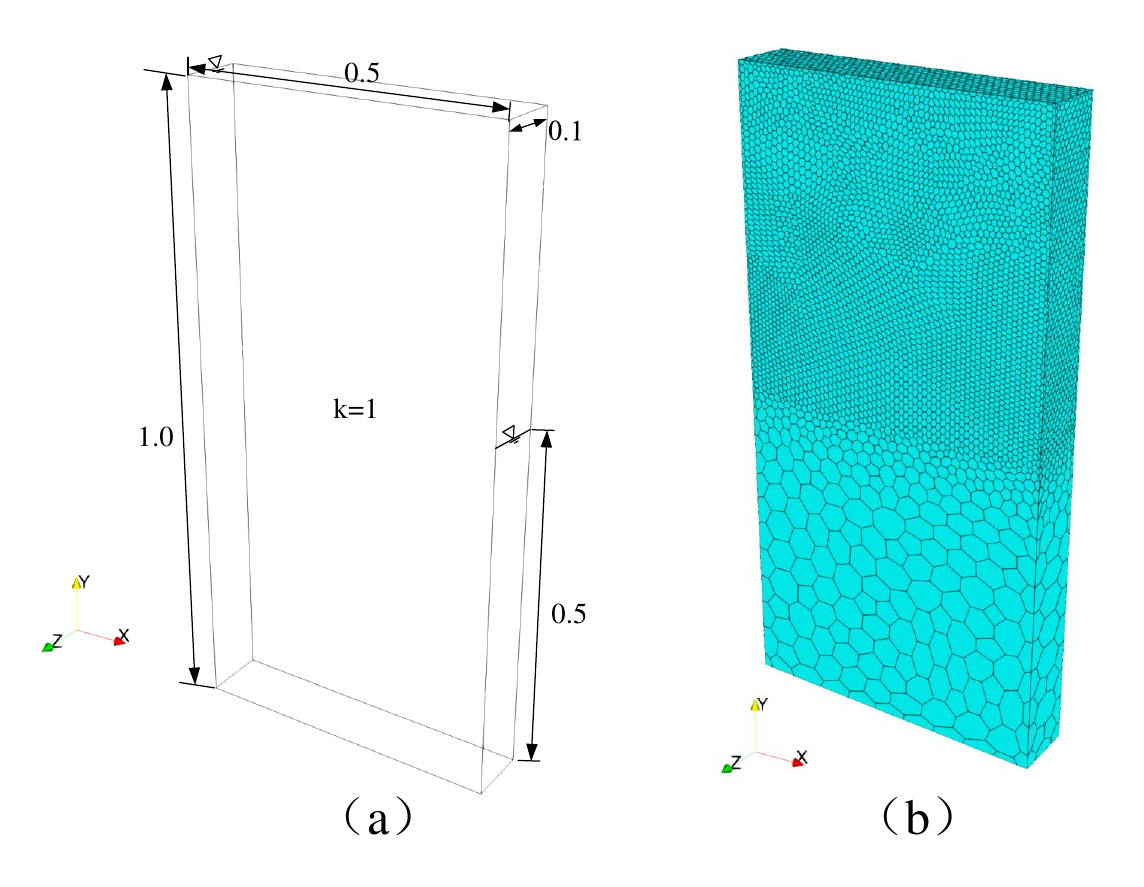}
   \caption{The diagram of a homogeneous rectangular dam; (a) geometry and boundary conditions; (b) mesh model.}
  \label{fig:ex05_geo_mesh}
\end{figure}

This example considers a homogeneous rectangular dam, with geometric details illustrated in Fig.~\ref{fig:ex05_geo_mesh}(a). The dam has a height of 1.0~m and a base width of 0.5~m. A hydraulic head of 1.0~m is applied along the upstream  boundary, while a lower head of 0.5~m is imposed at the downstream boundary. The bottom surface is treated as an impermeable boundary, and a constant permeability coefficient of \(K = 1~\mathrm{m/s}\) is assumed throughout the domain. To enhance the resolution near the expected location of the free surface above the downstream water level, a locally refined mesh is employed in this region, as depicted in Fig.~\ref{fig:ex05_geo_mesh}(b).

\begin{figure}[H]
  \centering
  \includegraphics[width=0.6\textwidth]{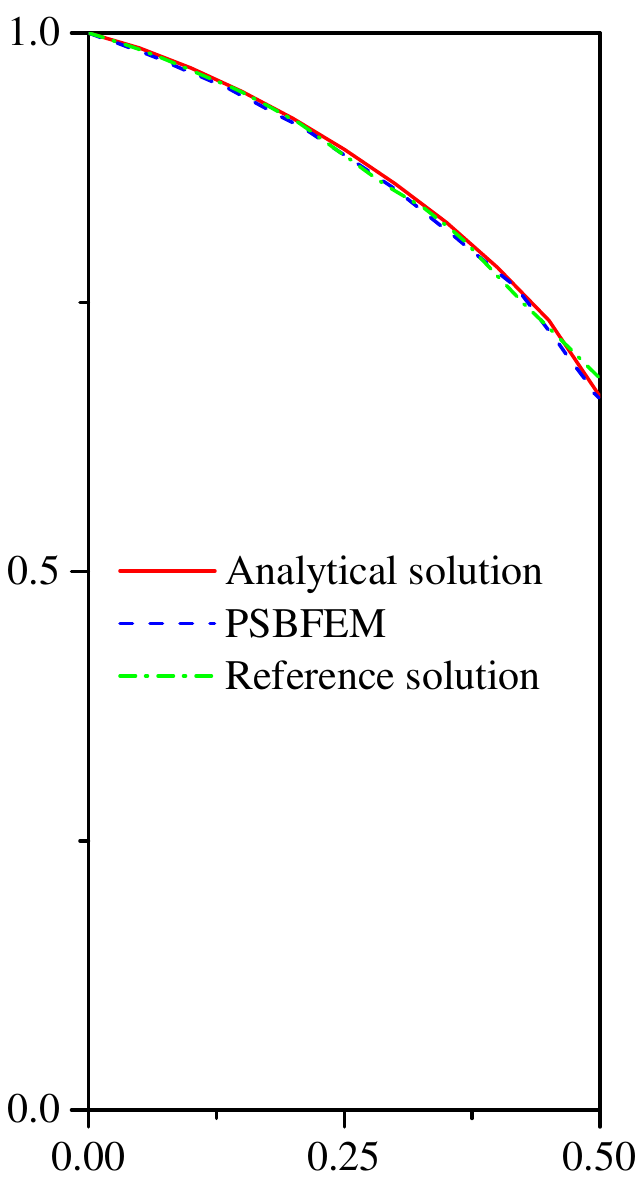}
   \caption{Comparison of free surface positions.}
  \label{fig:ex05_freeSurface}
\end{figure}

\begin{figure}[H]
  \centering
  \includegraphics[width=0.6\textwidth]{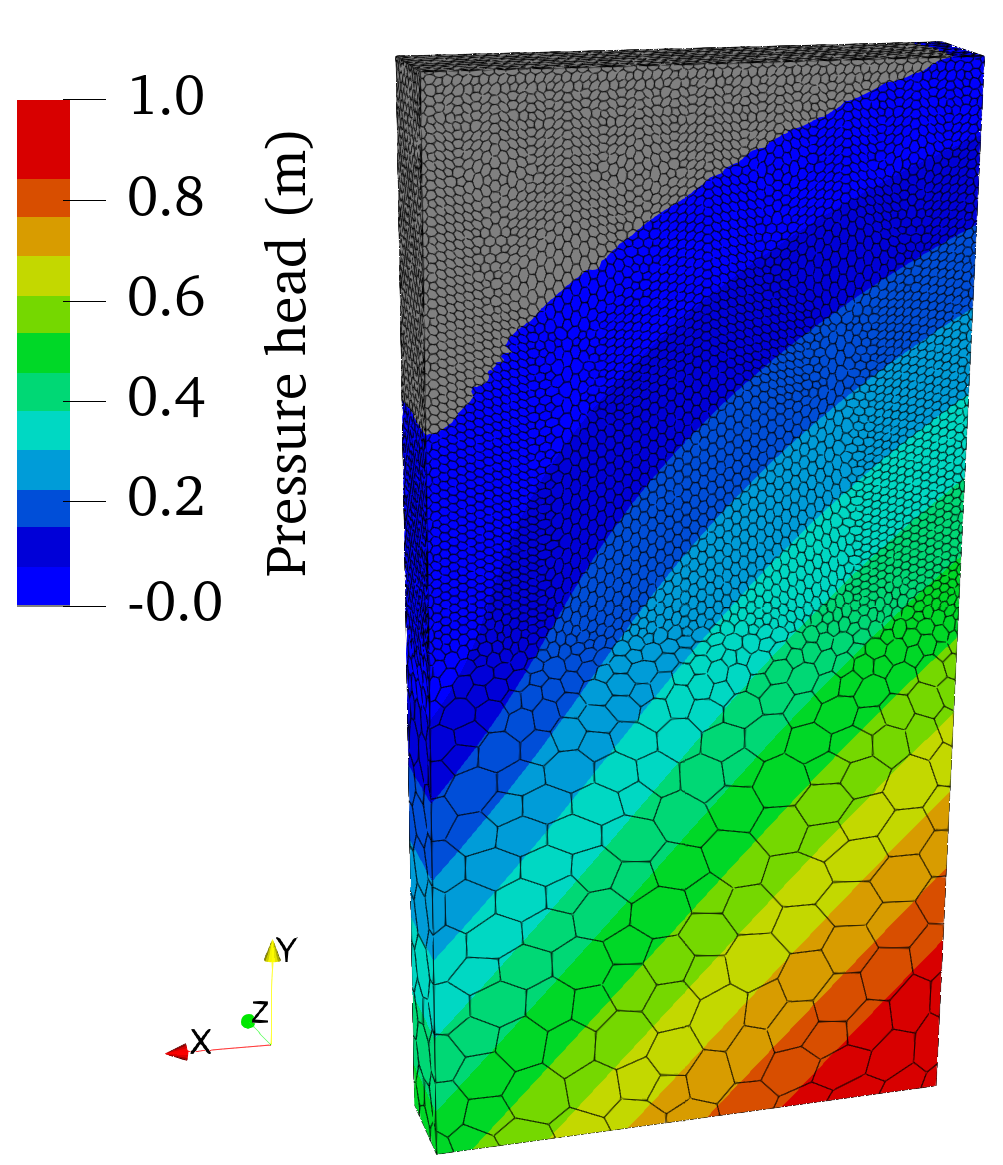}
   \caption{Distribution of pressure head of homogeneous rectangular dam.}
  \label{ex05_result_contour}
\end{figure}

Fig.~\ref{fig:ex05_freeSurface} presents the predicted free surface profiles obtained from various methods. The results demonstrate that the proposed PSBFEM solution closely aligns with the analytical reference. In addition, Fig.~\ref{ex05_result_contour} illustrates the computed pressure head distribution within the dam, further validating the capability of the PSBFEM to accurately capture seepage behavior in free-surface flow problems.

\begin{table}[htbp]
\centering
\caption{Coordinates of the overflow point.}
\begin{tabular}{lcc}
\toprule
Method & Coordinate of X (m) & Relative error\\
\midrule
Analytical solution & 0.662382 & - \\
PSBFEM & 0.661517 &$1.306\times10^{-3}$ \\
Reference solution \cite{jia2024new} & 0.671925&$1.441\times10^{-2}$ \\
\bottomrule
\end{tabular}
\label{tab:ex05}
\end{table}

\subsubsection{Homogeneous trapezoidal dam}
\begin{figure}[H]
  \centering
  \includegraphics[width=1.0\textwidth]{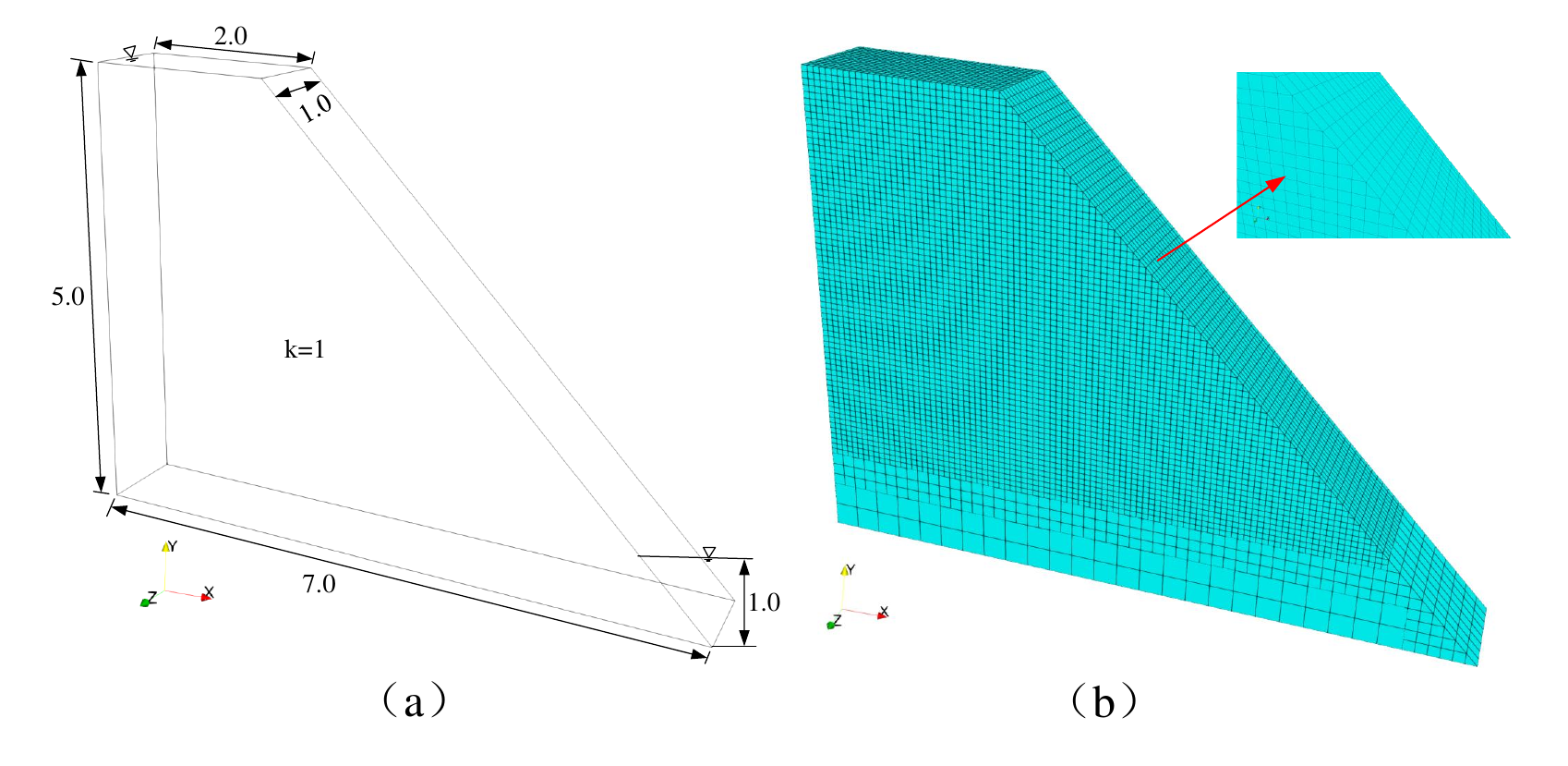}
   \caption{The diagram of a homogeneous trapezoidal dam; (a) geometry and boundary conditions; (b) mesh model.}
  \label{fig:ex06_geo_mesh}
\end{figure}

This example represents the cross-section of a homogeneous trapezoidal dam, with the computational model dimensions shown in Fig.~\ref{fig:ex06_geo_mesh} (a). The hydraulic head is 5 m on the left upstream side and 1 m on the right downstream side. The bottom is impermeable, and the permeability coefficient is 1 m/s. It can be observed from the  Fig.~\ref{fig:ex06_geo_mesh} (b) that the dam body is discretized using a hybrid octree mesh, with local refinement applied in the potential overflow region. It can be seen from the table that 98.67\% of the elements are octree elements, while only 1.53\% of the elements on the downstream slope are hybrid elements.

\begin{figure}[H]
  \centering
  \includegraphics[width=1.0\textwidth]{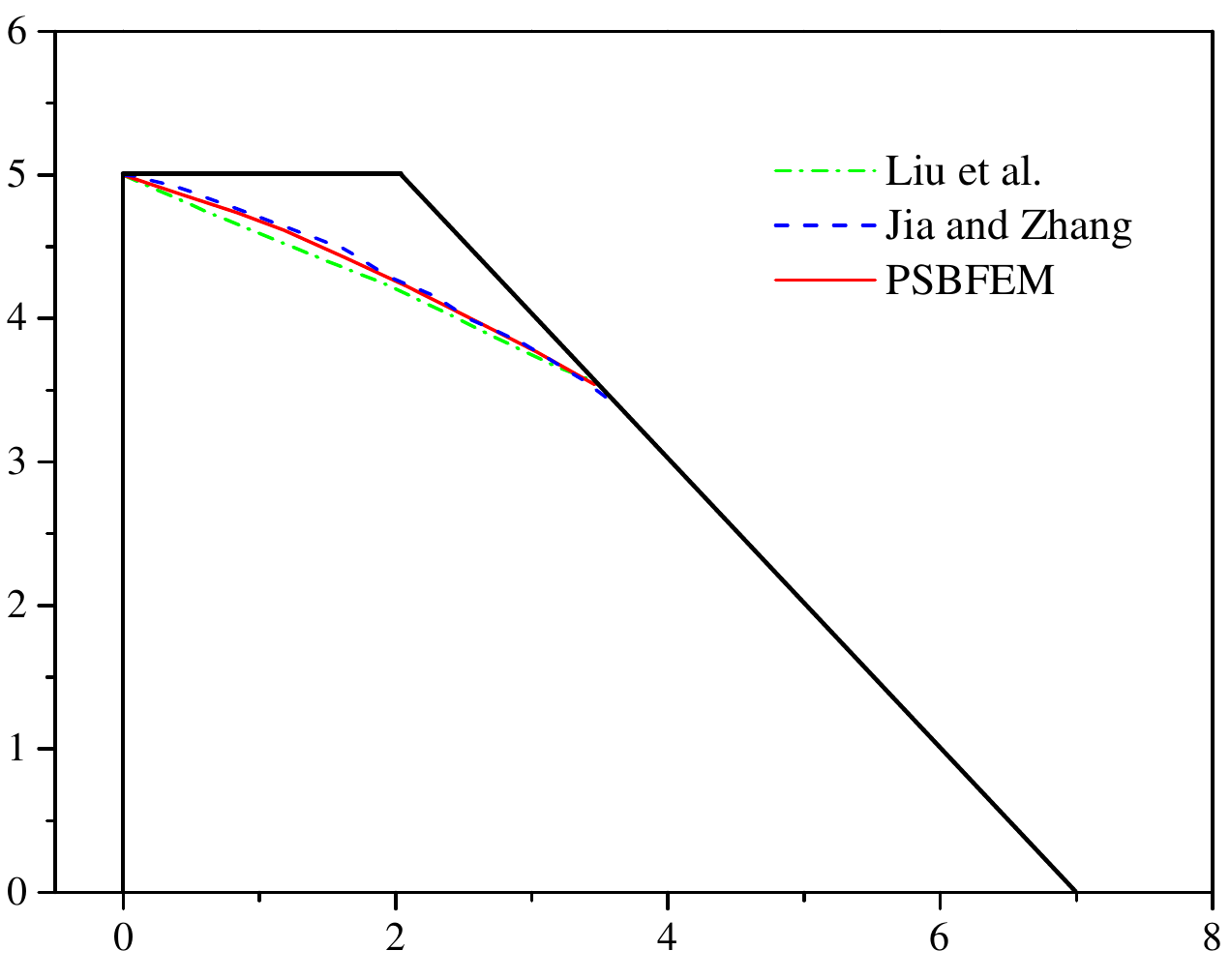}
   \caption{Comparison of free surface positions.}
  \label{fig:ex06_freeSurface}
\end{figure}

\begin{figure}[H]
  \centering
  \includegraphics[width=1.0\textwidth]{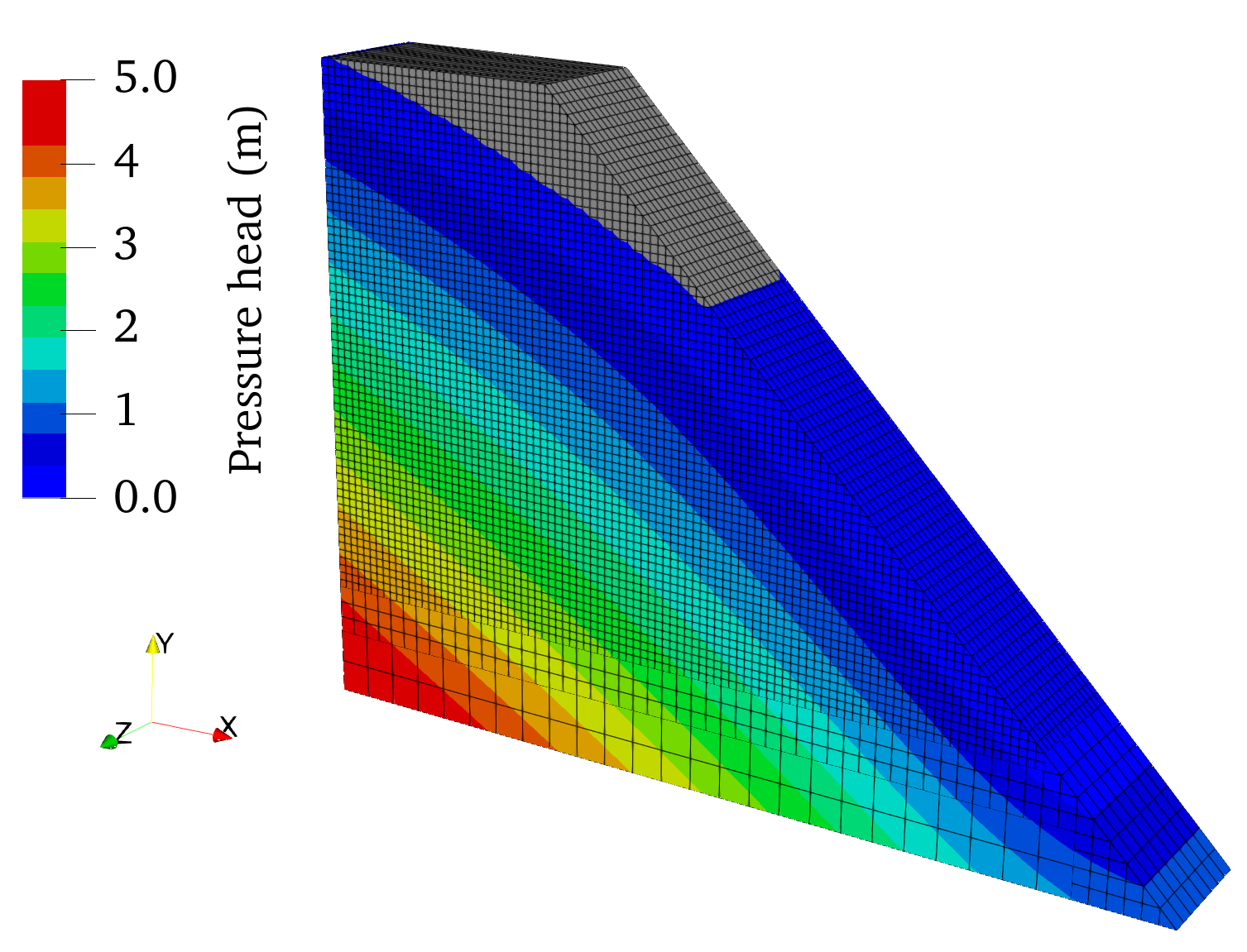}
   \caption{Distribution of pressure head for homogeneous trapezoidal dam.}
  \label{fig:ex06_result_contour}
\end{figure}

Fig.~\ref{fig:ex06_freeSurface} shows a comparison of free surface profiles predicted between the PSBFEM and FEM. The PSBFEM exhibits excellent consistency with the Liu et al.\cite{liu2018new} and Jia and Zheng \cite{jia2024new}. Additionally, the pressure head contours presented in Fig.~\ref{fig:ex06_result_contour} highlight the effectiveness of the method in capturing the seepage characteristics associated with free-surface flow.

\begin{table}[htbp]
\centering
\caption{Composition of hybrid octree mesh for the trapezoidal dam.}
\begin{tabular}{lccc}
\toprule
Element type & octree element & hybrid element\\
\midrule
Number of elements &69294	&1080 \\
Proportion &98.67\%	&1.53\% \\
\bottomrule
\end{tabular}
\label{tab:Composition of hybrid octree mesh}
\end{table}

\section{Conclusions}
\label{sec:8}
This work presented a PSBFEM framework for three dimensional seepage analysis. The developed method combines the semi-analytical formulation of the SBFEM with Wachspress shape functions for constructing conforming polygonal boundary interpolants, and employs polyhedral meshes for efficient volumetric discretization. This design yields significant advantages over conventional FEM in terms of geometric adaptability and numerical accuracy. The main conclusions are as follows:

(1) By employing polyhedral elements constructed via Wachspress interpolants, the PSBFEM improves mesh flexibility and enables accurate solutions with fewer elements. This substantially reduces preprocessing effort for complex geometries compared to conventional FEM, which often requires structured meshing or extensive manual adjustment.

(2) Numerical results show that the PSBFEM achieves high accuracy even on coarse or irregular meshes. Across steady-state, transient, and free-surface seepage scenarios, the proposed method consistently outperforms the FEM in convergence behavior, yielding lower relative errors under mesh refinement.

(3) The framework exhibits strong robustness when handling intricate domains, such as STL-based models. In particular, simulations involving hybrid octree meshes demonstrate that polyhedral discretization offers better accuracy and reduced computation time than FEM with tetrahedral meshes, underlining the suitability of PSBFEM for real-world geotechnical applications.

(4)  Notably, a fixed-mesh strategy is adopted for handling free-surface flows, which avoids remeshing while maintaining solution fidelity and numerical stability. In addition, the integration of hybrid octree meshes with SBFEM enables efficient local mesh refinement without compromising global accuracy. This approach offers a favorable balance between computational cost and precision, making it suitable for large-scale seepage problems with evolving boundaries.

We will focus on extending the PSBFEM framework to incorporate nonlinear seepage behaviors, unsaturated flow models, and multiphysics coupling such as hydro-mechanical interaction and thermo-seepage processes in future work. 

\section{Acknowledgements}
The National Natural Science Foundation of China (grant NO. 42167046), Yunnan Province Xing Dian Talent Support Program (grant NO. XDYC-QNRC-2022-0764), Yunan Funndamental Research Projects (grant NO. 202401CF070043) and the Science and Technology Talents and Platform Plan (grant NO. 202305AK34003) provided support for this study. 

 \bibliographystyle{elsarticle-num} 
 \bibliography{cas-refs}

\begin{thebibliography}{10}
\expandafter\ifx\csname url\endcsname\relax
  \def\url#1{\texttt{#1}}\fi
\expandafter\ifx\csname urlprefix\endcsname\relax\def\urlprefix{URL }\fi
\expandafter\ifx\csname href\endcsname\relax
  \def\href#1#2{#2} \def\path#1{#1}\fi

\bibitem{scudelerExam2017}
C.~Scudeler, C.~Paniconi, D.~Pasetto, M.~Putti, Examination of the seepage face boundary condition in subsurface and coupled surface/subsurface hydrological models 53~(3)  1799--1819.

\bibitem{zhouGroundwate2023}
C.-B. Zhou, Y.-F. Chen, R.~Hu, Z.~Yang, Groundwater flow through fractured rocks and seepage control in geotechnical engineering: {{Theories}} and practices 15~(1)  1--36.

\bibitem{shahSeepageLosses2021}
S.~A. Shah, M.~Kiran, T.~Khurshid, Seepage losses measurement of {{Desert}} minor and development of gauge-discharge rating curve: {{A}} case study in {{District Ghotki}}, {{Sindh}} 8~(1)  13--23.

\bibitem{johariStochas2019}
A.~Johari, A.~Talebi, Stochastic analysis of rainfall-induced slope instability and steady-state seepage flow using random finite-element method 19~(8)  04019085.

\bibitem{abokwiekFinite2022}
R.~Abokwiek, M.~Al~Sharabati, R.~Hawileh, J.~A. Abdalla, R.~Sabouni, G.~A. Husseini, A finite element model for the analysis of seepage flow of water under concrete dams 40~(5)  2823--2841.

\bibitem{liSelect2022}
X.~Li, X.~Li, Y.~Wu, L.~Wu, Z.~Yue, Selection criteria of mesh size and time step in {{FEM}} analysis of highly nonlinear unsaturated seepage process 146  104712.

\bibitem{schneiderLarge2022}
T.~Schneider, Y.~Hu, X.~Gao, J.~Dumas, D.~Zorin, D.~Panozzo, A large-scale comparison of tetrahedral and hexahedral elements for solving elliptic pdes with the finite element method, ACM Transactions on Graphics (TOG) 41~(3) (2022) 1--14.

\bibitem{zienkiewic2005}
O.~C. Zienkiewicz, R.~L. Taylor, J.~Z. Zhu, The Finite Element Method: Its Basis and Fundamentals, Elsevier.

\bibitem{hegen1996element}
D.~Hegen, Element-free galerkin methods in combination with finite element approaches, Computer Methods in Applied Mechanics and Engineering 135~(1-2) (1996) 143--166.

\bibitem{bonet2004variational}
J.~Bonet, S.~Kulasegaram, M.~Rodriguez-Paz, M.~Profit, Variational formulation for the smooth particle hydrodynamics (sph) simulation of fluid and solid problems, Computer Methods in Applied Mechanics and Engineering 193~(12-14) (2004) 1245--1256.

\bibitem{huang2020rkpm2d}
T.-H. Huang, H.~Wei, J.-S. Chen, M.~C. Hillman, Rkpm2d: an open-source implementation of nodally integrated reproducing kernel particle method for solving partial differential equations, Computational particle mechanics 7~(2) (2020) 393--433.

\bibitem{ghoneim2020smoothed}
A.~Y. Ghoneim, A smoothed particle hydrodynamics-phase field method with radial basis functions and moving least squares for meshfree simulation of dendritic solidification, Applied Mathematical Modelling 77 (2020) 1704--1741.

\bibitem{shojaei2022hybrid}
A.~Shojaei, A.~Hermann, C.~J. Cyron, P.~Seleson, S.~A. Silling, A hybrid meshfree discretization to improve the numerical performance of peridynamic models, Computer methods in applied mechanics and engineering 391 (2022) 114544.

\bibitem{luo2023novel}
T.~Luo, Y.~Feng, Q.~Huang, Z.~Zhang, M.~Yan, Z.~Yang, D.~Zheng, Y.~Yang, A novel solution for seepage problems using physics-informed neural networks, arXiv preprint arXiv:2310.17331 (2023).

\bibitem{cottrell2009isogeometric}
J.~A. Cottrell, T.~J. Hughes, Y.~Bazilevs, Isogeometric analysis: toward integration of CAD and FEA, John Wiley \& Sons, 2009.

\bibitem{gupta2023insight}
V.~Gupta, A.~Jameel, S.~K. Verma, S.~Anand, Y.~Anand, An insight on nurbs based isogeometric analysis, its current status and involvement in mechanical applications, Archives of Computational Methods in Engineering 30~(2) (2023) 1187--1230.

\bibitem{zeng2018smoothed}
W.~Zeng, G.~Liu, Smoothed finite element methods (s-fem): an overview and recent developments, Archives of Computational Methods in Engineering 25~(2) (2018) 397--435.

\bibitem{yan2024fast}
M.~Yan, Y.~Yang, C.~Su, Z.~Zhang, Q.~Duan, J.~Luo, G.~Xiong, T.~Luo, A fast cell-smoothed finite element method for solving static--dynamic problems using a hybrid quadtree mesh, International Journal of Computational Methods (2024) 2450071.

\bibitem{wu2023polygonal}
C.-T. Wu, S.-W. Wu, R.-P. Niu, C.~Jiang, G.~Liu, The polygonal finite element method for solving heat conduction problems, Engineering Analysis with Boundary Elements 155 (2023) 935--947.

\bibitem{yang2025steady}
Y.~Yang, M.~Yan, Z.~Zhang, D.~Hao, X.~Chen, W.~Chen, \href{https://arxiv.org/abs/2501.03908}{Steady-state and transient thermal stress analysis using a polygonal finite element method} (2025).
\newblock \href {http://arxiv.org/abs/2501.03908} {\path{arXiv:2501.03908}}.
\newline\urlprefix\url{https://arxiv.org/abs/2501.03908}

\bibitem{yan2025polyhedral}
M.~Yan, Y.~Yang, C.~Su, Z.~Zhang, Q.~Duan, D.~Hao, J.~Zhou, A polyhedral scaled boundary finite element method solving three-dimensional heat conduction problems, Engineering Analysis with Boundary Elements 175 (2025) 106191.

\bibitem{song2018scaled}
C.~Song, The scaled boundary finite element method: introduction to theory and implementation, John Wiley \& Sons, 2018.

\bibitem{krome2017semi}
F.~Krome, H.~Gravenkamp, A semi-analytical curved element for linear elasticity based on the scaled boundary finite element method, International Journal for Numerical Methods in Engineering 109~(6) (2017) 790--808.

\bibitem{yangDevelopmentABAQUSUEL2020}
Z.~Yang, F.~Yao, Y.~Huang, Development of {{ABAQUS UEL}}/{{VUEL}} subroutines for scaled boundary finite element method for general static and dynamic stress analyses, Engineering Analysis with Boundary Elements 114 (2020) 58--73.
\newblock \href {https://doi.org/10.1016/j.enganabound.2020.02.004} {\path{doi:10.1016/j.enganabound.2020.02.004}}.

\bibitem{zhang2024adaptive}
P.~Zhang, C.~Du, W.~Zhao, S.~Jiang, N.~Gong, N.~Bourahla, Z.~Qi, An adaptive sbfem based on a nonlocal macro/meso damage model for fracture simulation of quasibrittle materials, Engineering Fracture Mechanics 312 (2024) 110601.

\bibitem{ooi2010hybrid}
E.~Ooi, Z.~Yang, A hybrid finite element-scaled boundary finite element method for crack propagation modelling, Computer Methods in Applied Mechanics and Engineering 199~(17-20) (2010) 1178--1192.

\bibitem{li2016novel}
F.~Li, P.~Ren, A novel solution for heat conduction problems by extending scaled boundary finite element method, International Journal of Heat and Mass Transfer 95 (2016) 678--688.

\bibitem{chen2018efficient}
K.~Chen, D.~Zou, X.~Kong, X.~Yu, An efficient nonlinear octree sbfem and its application to complicated geotechnical structures, Computers and Geotechnics 96 (2018) 226--245.

\bibitem{zhang2020nonlocal}
Z.~Zhang, C.~Zhou, A.~Saputra, Z.~Yang, C.~Song, Nonlocal dynamic damage modelling of quasi-brittle composites using the scaled boundary finite element method, Engineering Fracture Mechanics 240 (2020) 107362.

\bibitem{gravenkamp2017efficient}
H.~Gravenkamp, A.~A. Saputra, C.~Song, C.~Birk, Efficient wave propagation simulation on quadtree meshes using sbfem with reduced modal basis, International Journal for Numerical Methods in Engineering 110~(12) (2017) 1119--1141.

\bibitem{zhang2024prismatic}
G.~Zhang, M.~Zhao, J.~Zhang, X.~Du, Prismatic-element sbpml coupled with sbfem for 3d infinite transient wave problems, Computer Methods in Applied Mechanics and Engineering 427 (2024) 117014.

\bibitem{ya2021open}
S.~Ya, S.~Eisentr{\"a}ger, C.~Song, J.~Li, An open-source abaqus implementation of the scaled boundary finite element method to study interfacial problems using polyhedral meshes, Computer Methods in Applied Mechanics and Engineering 381 (2021) 113766.

\bibitem{ye2021psbfem}
N.~Ye, C.~Su, Y.~Yang, Psbfem-abaqus: Development of user element subroutine (uel) for polygonal scaled boundary finite element method in abaqus, Mathematical Problems in Engineering 2021~(1) (2021) 6628837.

\bibitem{dai2015fully}
S.~Dai, C.~Augarde, C.~Du, D.~Chen, A fully automatic polygon scaled boundary finite element method for modelling crack propagation, Engineering Fracture Mechanics 133 (2015) 163--178.

\bibitem{natarajan2017scaled}
S.~Natarajan, E.~T. Ooi, A.~Saputra, C.~Song, A scaled boundary finite element formulation over arbitrary faceted star convex polyhedra, Engineering Analysis with Boundary Elements 80 (2017) 218--229.

\bibitem{saputra2020three}
A.~A. Saputra, S.~Eisentr{\"a}ger, H.~Gravenkamp, C.~Song, Three-dimensional image-based numerical homogenisation using octree meshes, Computers \& Structures 237 (2020) 106263.

\bibitem{saputra2017automatic}
A.~Saputra, H.~Talebi, D.~Tran, C.~Birk, C.~Song, Automatic image-based stress analysis by the scaled boundary finite element method, International Journal for Numerical Methods in Engineering 109~(5) (2017) 697--738.

\bibitem{li2012scaled}
F.~Li, Q.~Tu, The scaled boundary finite element analysis of seepage problems in multi-material regions, International Journal of Computational Methods 9~(01) (2012) 1240008.

\bibitem{liu2018new}
J.~Liu, J.~Li, P.~Li, G.~Lin, T.~Xu, L.~Chen, New application of the isogeometric boundary representations methodology with sbfem to seepage problems in complex domains, Computers \& Fluids 174 (2018) 241--255.

\bibitem{johari2018reliability}
A.~Johari, A.~Heydari, Reliability analysis of seepage using an applicable procedure based on stochastic scaled boundary finite element method, Engineering Analysis with Boundary Elements 94 (2018) 44--59.

\bibitem{yangNovelSolutionSeepage2022}
Y.~Yang, Z.~Zhang, Y.~Feng, K.~Wang, A {{Novel Solution}} for {{Seepage Problems Implemented}} in the {{Abaqus UEL Based}} on the {{Polygonal Scaled Boundary Finite Element Method}}, Geofluids 2022 (2022) 1--14.
\newblock \href {https://doi.org/10.1155/2022/5797014} {\path{doi:10.1155/2022/5797014}}.

\bibitem{yan2023psbfem}
M.~Yan, Y.~Yang, Z.~Zhang, C.~Su, T.~Luo, A psbfem approach for solving seepage problems based on the pixel quadtree mesh, Geofluids 2023~(1) (2023) 9092488.

\bibitem{songScaled1999}
C.~Song, J.~P. Wolf, The scaled boundary finite element method-alias consistent infinitesimal finite element cell method-for diffusion, International Journal for Numerical Methods in Engineering 45~(10) (1999) 1403--1431.
\newblock \href {https://doi.org/10.1002/(SICI)1097-0207(19990810)45:10<1403::AID-NME636>3.0.CO;2-E} {\path{doi:10.1002/(SICI)1097-0207(19990810)45:10<1403::AID-NME636>3.0.CO;2-E}}.

\bibitem{yangDevelopment2020}
Z.~Yang, F.~Yao, Y.~Huang, Development of {{ABAQUS UEL}}/{{VUEL}} subroutines for scaled boundary finite element method for general static and dynamic stress analyses, Engineering Analysis with Boundary Elements 114 (2020) 58--73.
\newblock \href {https://doi.org/10.1016/j.enganabound.2020.02.004} {\path{doi:10.1016/j.enganabound.2020.02.004}}.

\bibitem{wachspress2006rational}
E.~L. Wachspress, A rational basis for function approximation, in: Conference on Applications of Numerical Analysis: Held in Dundee/Scotland, March 23--26, 1971, Springer, 2006, pp. 223--252.

\bibitem{warren2003uniqueness}
J.~Warren, On the uniqueness of barycentric coordinates, Contemporary Mathematics 334 (2003) 93--100.

\bibitem{warren2007barycentric}
J.~Warren, S.~Schaefer, A.~N. Hirani, M.~Desbrun, Barycentric coordinates for convex sets, Advances in computational mathematics 27 (2007) 319--338.

\bibitem{yu2021scaled}
B.~Yu, P.~Hu, A.~A. Saputra, Y.~Gu, The scaled boundary finite element method based on the hybrid quadtree mesh for solving transient heat conduction problems, Applied Mathematical Modelling 89 (2021) 541--571.

\bibitem{zienkiewicz1989finite}
O.~Zienkiewicz, R.~Taylor, O.~Zienkiewicz, R.~Taylor, The finite element method, vol. 1 mcgraw-hill (1989).

\bibitem{ANSYSFluent}
{ANSYS Inc.}, ANSYS Fluent User’s Guide, release 2023 R1, Canonsburg, PA, USA (2023).

\bibitem{STARCCM}
{Siemens Digital Industries Software}, STAR-CCM+ User Guide, version 2023.1, Siemens PLM Software Inc., Plano, TX, USA (2023).

\bibitem{ParaViewGuide}
{Kitware, Inc.}, The ParaView Guide: A Parallel Visualization Application, version 5.11, Clifton Park, NY, USA (2023).

\bibitem{bathe1979finite}
K.-J. Bathe, M.~R. Khoshgoftaar, Finite element free surface seepage analysis without mesh iteration, International Journal for Numerical and Analytical Methods in Geomechanics 3~(1) (1979) 13--22.

\bibitem{zhang2019fast}
J.~Zhang, S.~Chauhan, Fast explicit dynamics finite element algorithm for transient heat transfer, International Journal of Thermal Sciences 139 (2019) 160--175.

\bibitem{yangNovel2022}
Y.~Yang, Z.~Zhang, Y.~Feng, K.~Wang, A {{Novel Solution}} for {{Seepage Problems Implemented}} in the {{Abaqus UEL Based}} on the {{Polygonal Scaled Boundary Finite Element Method}}, Geofluids 2022 (2022) 1--14.
\newblock \href {https://doi.org/10.1155/2022/5797014} {\path{doi:10.1155/2022/5797014}}.

\bibitem{venus_of_arles_3axis}
{3Axis.co}, \href{https://3axis.co/venus-of-arles-head-3d-printer-model/qwkdg5zqkl/}{Venus of arles head - 3d printer model}, accessed: 2025-03-13 (2025).
\newline\urlprefix\url{https://3axis.co/venus-of-arles-head-3d-printer-model/qwkdg5zqkl/}

\bibitem{jia2024new}
Z.~Jia, H.~Zheng, A new procedure for locating free surfaces of complex unconfined seepage problems using fixed meshes, Computers and Geotechnics 166 (2024) 106032.

\end{thebibliography}





\end{document}